\documentclass{article}
\usepackage{amsmath,amssymb,latexsym,enumerate}
\newtheorem{lemma}{Lemma}
\newtheorem{theorem}{Theorem}
\newtheorem{corollary}{Corollary}

\renewcommand{\Pr}{{\rm Pr}}

\newcommand{\E}{{\rm E}}
\newcommand{\M}{\mathcal{M}}

\newcommand{\T}{\mathcal{T}}

\newcommand{\old}[1]{{}}

\begin{document}
\title{A new asymptotic enumeration technique: the  Lov\'asz Local Lemma }
\author{Linyuan Lu \thanks{
This researcher was supported in part by the NSF DMS contracts Nos. 0701111, 1000475  and 1300547.}
\ \ \ \ 
L\'aszl\'o Sz\'ekely \thanks{
This researcher was supported in part by
the NSF DMS contracts  Nos. 0701111, 1000475,  and 1300547, by the Alexander von Humboldt Foundation at the Rheinische
     Friedrich-Wilhelms Universit\"at, Bonn, and the contract \#FA9550-12-1-0405 from the U.S.
Air Force Office of Scientific Research (AFOSR) and the Defense Advanced
Research Projects Agency (DARPA).}\\
University of South Carolina {\tt \{lu, szekely\}@math.sc.edu}\\
}

\maketitle

\begin{abstract}
Our previous paper \cite{soda} applied a lopsided version of the Lov\'asz Local
Lemma that allows negative dependency graphs \cite{erdosspencer} to the space
of random injections from an $m$-element set to an $n$-element set.
(Equivalently, the same story can be told about the space of random matchings
in $K_{n,m}$.) In this paper we show how the lopsided version of the  Lov\'asz Local
Lemma applies to the space of random matchings
in $K_{2n}$. We also prove  tight upper bounds
that  asymptotically match the lower bound given by the
 Lov\'asz Local Lemma. As a consequence, we give new proofs to a number of results
on the enumeration of permutations, Latin rectangles, and regular graphs.
The strength of the method is shown by a new result:
enumeration of graphs by degree sequence or bipartite degree sequence and girth. 
As another application, we provide a new proof to the classical
probabilistic result of Erd\H os \cite{erdoscan}
 that showed the existence of graphs 
with arbitrary large girth and chromatic number. 
If the degree sequence satisfies some mild conditions, almost all graphs with this
degree sequence and prescribed girth have high chromatic number. 
\end{abstract}

\section{Lov\'asz Local Lemma with negative dependency graphs}
This is a sequel to our previous paper \cite{soda} and 
we use the same notations.
Let $A_1,A_2,\ldots, A_n$ be events in a probability space.

A {\em negative dependency graph} for $A_1,\ldots, A_n$
is a simple graph on $[n]$ satisfying
\begin{equation}
   \label{eq:neg2}\Pr(A_i|\wedge_{j\in S} \overline{ A_j})\leq \Pr(A_i),
\end{equation}
for any index $i$ and any subset $S\subseteq \{j\mid ij \not \in
E(G)\}$,  whenever the conditional probability 
 $\Pr(A_i \mid  \wedge_{j\in S} \overline{ A_j})$ is well-defined, i.e.
 $\Pr( \wedge_{j\in S} \overline{ A_j})>0$.
We will make use of the fact that
inequality (\ref{eq:neg2}) trivially holds when $\Pr(A_i)=0$, otherwise
the following inequality is equivalent to inequality
(\ref{eq:neg2}):
\begin{equation}
   \label{eq:neg1}
   \Pr(\wedge_{j\in S} \overline{A_j} \mid A_i)\leq \Pr(\wedge_{j\in S}
\overline{ A_j}).
\end{equation}

\noindent
For variants of the Lov\'asz Local Lemma with increasing strength,
see \cite{erdoslovasz, spencerbook, erdosspencer,ku}:
\begin{lemma}
{\bf [Lov\'asz Local Lemma.]}  
\label{third}
Let $A_1,\ldots, A_n$ be events with a negative dependency graph $G$.
If there exist numbers $x_1,\dots, x_n\in [0,1)$ such that
\begin{equation}
  \label{eq:ai}
\Pr(A_i)\leq x_i\prod_{ij\in E(G)}(1-x_j)  
\end{equation}
for all $i$, then
\begin{equation} \label{kovetk}
\Pr(\wedge_{i=1}^n \overline{A_i})\geq \prod_{i=1}^n(1-x_i).
\end{equation}
\end{lemma}

The main obstacle for using Lemma \ref{third} is the difficulty to
define a useful negative dependency graph other than a dependency
graph.
In \cite{soda}, we described a general way to create negative dependency
graphs in the space of random functions $U\rightarrow V$ equipped with
uniform distribution. Namely, let the events be the set of 
all extensions of some particular
{\em partial} functions to functions; 
and create an edge for the negative dependency graph,
if the partial functions have common elements in their domains or ranges,
other than the agreement of the partial functions.
These events also can be thought of as all extensions of (partial) matchings 
in the complete bipartite graph with classes $U,V$, where an edge of the
negative dependency graph comes from two event-defining 
(partial) matchings whose
union is no longer a (partial) matching. 
 In \cite{soda}, we used this technique
to prove a new result on the Tur\'an hypergraph problem, and we found
surprising applications as proving lower bounds (matching 
certain   asymptotic formulas) for permutation and Latin rectangle enumeration
problems.

In this paper, we show an analogous construction of a  
negative dependency graph for events, which live in the space of random 
matchings of a complete graph. We require that the events are the set of
all extensions of (partial) matchings
in  a complete graph to perfect matchings, and two event-defining
partial matchings make an edge, if their union is no longer a (partial)
matching. 
 Our construction, however, fails to provide negative dependency graphs for extensions of
partial matchings of arbitrary graphs.

We move one step further and show some general and some specific {\em upper
bounds} for the event estimated by the Lov\'asz Local Lemma, and show that
for large classes of problems the upper bound is asymptotically equal to
the lower bound. These results apply to the permutation enumeration problems
in \cite{soda}, and to enumeration problems of regular graphs.
Many
asymptotic enumeration results
that we prove are not new and typically do not give
 the largest known valid
range of the asymptotic formula, but are nontrivial results and
 often more recent
than the Lov\'asz Local Lemma itself. They come
out from our framework elementarily, and even easily. 
 
 The strength of the framework is shown by a new result:
enumeration of graphs by degree sequence and girth, under mild condititions 
for the degree sequence. We also provide an analogous enumeration result
for bipartite degree sequence and girth. Although they are special cases, we 
prove the results for regular graphs first, as they simplify the explanation
for more general degree sequences.
 
 There is literature on some  improvements on the  Lov\'asz Local Lemma
 using methods of statistical physics, e.g. \cite{scott}, \cite{procacci}, that we do 
 not touch upon this paper, as they are difficult to use and the improvement would be tiny,
 if present at all, in a resulting  asymptotic formula.
 

As another application, we revisit a classic of the probabilistic method:
Erd\H os' proof to the existence of graphs with arbitrary large girth and 
chromatic number \cite{erdoscan}.  We show that if
 the degree sequence satisfies some mild conditions, almost all graphs with this
degree sequence and prescribed girth have high chromatic number. 

In a scenario of the  Poisson paradigm, we  estimate the probability that
none of a set of rare events occur. Let $X$ be the sum of  the indicator
variables 
of these events and $\mu=\E(X)$. If the dependency among these events
is rare, then one would expect that $X$ has a Poisson distribution
with mean $\mu$. In particular, $\Pr(X=0)\approx e^{-\mu}$.
The Janson inequality and Brun's sieve method \cite{as}
are often the good choice to solve these kind of  problems.
Now we offer an alternative approach---using
Lov\'asz Local Lemma. Our approach can be considered as an analogue
of the Janson inequality in another setting that offers plenty of 
applications. It is curious that the proof of Boppana and Spencer
\cite{boppanas} for the
Janson inequality (see also in \cite{as}) uses conditional probabilities 
somewhat similarly to the proof of the Lov\'asz Local Lemma. There is an inherent 
reason why we do not get the "second term" in asymptotic enumeration,
like in (\ref{latinas}) or (\ref{strongest}), which extends the range of the asymptotic formula: $e^{-\mu}$
is {\em between} our lower and upper bounds (see Theorem~\ref{enum}), and therefore
we cannot add a correction term to $-\mu$ in the exponent.

For further research, it would be interesting to get asymptotics for
further terms from the Poisson distribution, i.e. for the probability of
exactly $k$ events holding, for any fixed $k$. 
Lots  of further applications of our framework are possible, this paper 
gives just a sampler of applications.

\section{Some general results on negative and near-positive 
dependency graphs} \label{genbounds}
These lower and upper bounds are {\em general} in the sense that there is no
assumption on the events being defined through matchings.

All over this paper, we will be using properties of a useful function, which cannot be expressed 
in terms of elementary functions, but can be expressed with LambertW. Recall that LambertW is a multivalued function 
satisfying $$z=\hbox{LambertW}(z)e^{\hbox{LambertW}(z)}.$$
In the following lemma we summarize the properties that we will need.
\begin{lemma} \label{function}
\begin{enumerate}[{\rm (i)}]
\item For $0\leq \gamma \leq 1/4$, the equation 
\begin{equation} \label{defeq}
1=ye^{-\gamma y}
\end{equation}
has a unique solution  $y$ in $1\leq y\leq 2$, and defines a function $y(\gamma)$.
\item $y(\gamma)=-\frac{{W_0}(-\gamma)}{\gamma}$, where  $W_0$ is the branch of {\rm LambertW} with $W_0(0)=0$.
\item As  the Taylor series of ${W_0}(\gamma)$  around 0 is convergent for $|\gamma|<1/e$,
so is    the Taylor series of $y(\gamma)$ around 0. 
\item $y(\gamma)$ is strictly increasing on $[0, 1/4]$.
\item For $\gamma\rightarrow 0$,
\begin{equation} \label{lagbur}
y(\gamma)=1+\gamma +\frac{3}{2}\gamma^2+\frac{8}{3}\gamma^3+\frac{125}{24}\gamma^4
+\frac{54}{5}\gamma^5+O(\gamma^6).
\end{equation}
\item For $0\leq \gamma \leq 1/4$, 
\begin{equation} \label{lagburexact}
 1+\gamma +\frac{3}{2}\gamma^2\leq y(\gamma)\leq 1+\gamma +\frac{3}{2}\gamma^2+66\gamma^3.
\end{equation}
\end{enumerate}
\end{lemma}
{\bf Proof:} (ii) and (v) can be obtained with Maple. 
As the RHS of (\ref{defeq}) $<1$ at $y=1$ and $>1$ at $y=2$, there is a solution in between
for (\ref{defeq}), providing the existence for (i).
Using implicit differentiation, $y'(\gamma) >0 $ in $[0, 1/4]$, proving (iv) and the uniqueness
claim in (i). Finally, for (vi), estimates for $y'''(\gamma)$ were obtained with Maple.
\hfill$\square$

Many results in this paper are of asymptotic nature. Assume that for all (or infinitely many) positive integers $N$ 
there is a probability space  \linebreak  $(\Omega(N),{\cal A}(N), \Pr_N)$ and events $A_1(N),...,A_{n(N)}(N)\in  {\cal A}(N)$. We 
consider a sequence of problems: 
obtain estimates or asymptotic formula for  $$\Pr_N\bigl(\wedge_{i=1}^{n(N)} \overline{A_i(N)}\bigl).$$
The use of little-oh or  big-Oh formulae  and asymptotics all refer to $N\rightarrow \infty$. For simplicity, however,
from now on we do not make  $N$ explicit in the notation.
In many sequences of problems  
  $\Pr(A_i)$ and $\sum_{ij\in E(G)} \Pr(A_j)$
are  so small  that one can set
$x_i=:(1+o(1))\Pr(A_i)$ to use Lemma \ref{third}. 
\begin{theorem}\label{simple}
Let $A_1,\ldots, A_n$ be events with negative dependency graph $G$.
Let us be given any $\epsilon$ with $0<\epsilon<1/4$.
If  
\begin{equation}  \label{altcond}\Pr(A_i)<\epsilon \hbox{\ \  and \ \ } 
\sum_{j:ij\in E(G)} \Pr(A_j)+ 2\Pr^2(A_j)<\epsilon
\end{equation}
  for every  $1\leq i\leq n$,  then
\begin{enumerate}[{\rm (i)}]
\item 
for any $S,T\subseteq V(G)$ with $S\cap T=\emptyset$,
 we have
 \begin{equation}
\Pr(\wedge_{i\in S} \overline{A_i} \mid\wedge_{j\in T}  \overline{A_j}) \geq 
\prod_{i\in S} \Bigl(1- \Pr(A_i) y(\epsilon ) \Bigl);
\label{l5formula}
 \end{equation}
\item in particular, we have  
\begin{equation} \label{simple1}
 \Pr(\wedge_{i=1}^n \overline{A_i})\geq 
\exp\Biggl({-\sum_{i=1}^n\Pr(A_i) y(\epsilon)-   \sum_{i=1}^n\Pr^2(A_i) y^2(\epsilon)}\Biggl) .
\end{equation}
 \end{enumerate}
\end{theorem}
{\bf Proof:} Set $x_i=\Pr(A_i) y(\epsilon).$ It is clear that
$0\leq x_i <1/2$. Observe that for $0\leq x\leq 1/2$ we have $1-x\geq e^{-x-x^2}$. To
use Lemma~\ref{third}, we need the condition (\ref{eq:ai}). Indeed, 
$\Pr(A_i)= x_i/y(\epsilon)=x_ie^{-\epsilon y(\epsilon)}  \leq x_i\exp\bigl(-\sum_{j:ij\in E(G)}
(x_j+x_j^2)\bigl)\leq x_i\prod_{j:ij\in E(G)}(1-x_j).$
To prove (i), we recall  not the conclusion of Lov\'asz Local Lemma, 
but a crucial step
in the proof (see \cite{spencerbook}, \cite{ku}): for any $T\subseteq V(G)$ with 
$i\notin T$, we have 
$\Pr(A_i \mid \wedge_{j\in \T, j\not= i} \overline{A_j})\leq x_i$,
which in our case yields for any $i\in S$
\[\Pr(A_i\mid \wedge_{j'\in T} \overline{A_{j'}})\leq x_i
=\Pr(A_i )y(\epsilon)
.\]
Assume that $S=\{m_1,m_2,...,m_s\}$. We have 
\[
\Pr(\overline{A_{m_1}}\wedge \overline{A_{m_2}} \wedge .... \wedge 
\overline{A_{m_s}}\mid\wedge_{j\in T}  \overline{A_j}) =
\]
\[
\prod_{\ell=1}^s \Biggl[\Pr\biggl(\overline{A_{m_\ell }} \mid
\overline{A_{m_1}}\wedge \overline{A_{m_2}} \wedge .... \wedge 
\overline{A_{m_{\ell-1 }}}\wedge (\wedge_{j\in T}  \overline{A_j})
\biggl)\Biggl]=
\]
\[
\prod_{\ell=1}^s \Biggl[1-\Pr\biggl(A_{m_\ell } \mid
\overline{A_{m_1}}\wedge \overline{A_{m_2}} \wedge .... \wedge 
\overline{A_{m_{\ell-1 }}}\wedge (\wedge_{j\in T}  \overline{A_j})
\biggl)\Biggl]\geq \prod_{\ell =1}^s (1-x_{m_\ell}).
\]
 The conclusion  of (ii) is implied by (i) with $T=\emptyset$ or by  Lemma~\ref{third}:  $\Pr(\wedge_{i=1}^n \overline{A_i})\geq \prod_i(1-x_i)=\prod_i\bigl(1-\Pr(A_i)y(\epsilon)\bigl)\geq \exp\Bigl({-\sum_{i=1}^n\Pr(A_i) y(\epsilon)-   \sum_{i=1}^n\Pr^2(A_i) y^2(\epsilon)}\Bigl).$
 \hfill$\square$

\noindent Theorem~\ref{simple} provided {\sl logarithmic asymptotics} for  the expected Poisson type lower
bound when $\epsilon \rightarrow 0$ for a {\sl sequence} of problems and estimations. However,
we want {\em asymptotics},     and obtain it with slightly more assumptions:
\begin{corollary}\label{LLLlower}
 Set $\mu=\sum_i \Pr(A_i)$. If  for a  sequence of problems  $\epsilon\mu\rightarrow 0$, then
 \begin{equation}
\Pr(\wedge_{i=1}^n \overline{A_i})\geq (1-o(1))e^{-\mu}.
 \end{equation}
 This holds, in particular, when $\mu$  is bounded and 
 $\epsilon\rightarrow 0$.
\end{corollary}

We comment here that this result does not allow a good generalization with {\sl different} 
bounds on $\Pr(A_i)$ and $\sum_{j:ij\in E(G)} \Pr(A_j)$.

Next we give a crucial new definition.
For the events  $A_1,\ldots, A_n$  in a probability space $\Omega$,
and an $\epsilon$ with
 $1>\epsilon>0$, we define an {\it  $\epsilon$-near-positive 
dependency graph}
to be a graph $G$  on $V(G)=[n]$ satisfying
\begin{enumerate}[{\rm (i)}]
\item $\Pr(A_i\wedge A_j)=0$ if $ij\in E(G)$.
\item For any index $i$ and any subset $i\notin T\subseteq \{j\mid ij \not \in
E(G)\}$,
\[\Pr(A_i\mid \wedge_{j\in T} \overline{ A_j} ) \geq (1-\epsilon)\Pr(A_i),\]
whenever the conditional probability is well-defined.
\end{enumerate}

\begin{theorem}\label{ub}
 Let $A_1,\ldots, A_n$ be events with an $\epsilon$-near-positive
dependency graph $G$.
Then we have
\[\Pr(\wedge_{i=1}^n \overline{A_i})\leq
\prod_{i=1}^n [1-(1-\epsilon)\Pr(A_i)].\]
\end{theorem}
{\bf Proof:} If $\Pr(\wedge_{i=1}^n \overline{A_i})=0$, then the
conclusion holds. So we may assume without loss of generality that
$\Pr(\wedge_{i=1}^n \overline{A_i})>0$.
Now we would like to show that for any $i$ and any subset $S\subseteq V(G)$
 with $i\notin S$,
\[\Pr(A_i\mid \wedge_{j\in S} \overline{ A_j})\geq (1-\epsilon)\Pr(A_i),\]
as the conditional probability above is well-defined by our assumption.
Write $S=S_1\cup S_2$, where $S_1=S\cap N_G(i)$ and $S_2=S\setminus S_1$.
We have
\begin{eqnarray*}
  \Pr(A_i\mid \wedge_{j\in S} \overline{ A_j})
&=& \frac{\Pr(A_i\wedge(\wedge_{k\in S_1}\overline{
 A_k}) \mid \wedge_{j\in S_2} \overline{ A_j})}
{\Pr(\wedge_{k\in S_1}\overline{ A_k}  \mid \wedge_{j\in S_2} \overline{ A_j})}\\
&=& \frac{\Pr(A_i\mid \wedge_{j\in S_2} \overline{ A_j})}
{\Pr(\wedge_{k\in S_1}\overline{ A_k}  \mid \wedge_{j\in S_2} \overline{ A_j}
)}\\
&\geq&
\Pr(A_i\mid \wedge_{j\in S_2} \overline{ A_j}) \\
&\geq& (1-\epsilon)\Pr(A_i).
\end{eqnarray*}
(The first part of the definition of the $\epsilon$-near-positive 
dependency graph, $\Pr(A_i\wedge A_j)=0$ for $ij$ edges,
allowed the elimination of the  $\wedge_{k\in S_1}\overline{A_k} $ term.)
Hence, we have
\begin{eqnarray*}
  \Pr(\wedge_{i=1}^n \overline{A_i})
&=&\prod_{i=1}^n \Pr(\overline{ A_i} \mid \wedge_{k=i+1}^n \overline{ A_k})=\\
\prod_{i=1}^n [1-\Pr({ A_i} \mid \wedge_{k=i+1}^n \overline{ A_k})]
&\leq& \prod_{i=1}^n (1-(1-\epsilon)\Pr(A_i)). \hbox{\ \ \ \ \ \ \ \ \ }\square
\end{eqnarray*}

\section{Instances for negative dependency graphs:
The space of random matchings of $K_N$ and $K_{N,N'}$}
Let $\Omega$ denote the probability space of
perfect matchings  of the complete bipartite graph $K_{N,N'}$ ($N\leq N'$)
or the  probability space of  the complete graph $K_{N}$  for an even integer $N$;
 equipped with the
uniform distribution.  We are going to apply the  Lov\'asz Local Lemma
(Lemma \ref{third})
in $\Omega$ by identifying a class of negative dependency graphs.  For any
(not necessary perfect) matching $M$, let $A_M$ be the set
of   maximum cardinality (in $K_N$ perfect) matchings extending $M$:
\begin{equation}
  \label{eq:1}
  A_M=\{F\in \Omega \mid M\subseteq F\}.
\end{equation}
We will term an event $A_M$  in (\ref{eq:1}), with $M\not= \emptyset$, 
a {\it canonical event}. 
We will say that two matchings, $M_1$ and $M_2$,
 are in {\it conflict}, if $M_1\cup M_2$ is not a matching.
For a matching $M$, we will denote by $supp(M)$ the support set of
the matching, i.e. the
$2|M|$ vertices that its edges cover. We leave the following easy lemma to the reader:
\begin{lemma} \label{onconflict}
\begin{enumerate}[{\rm (i)}]
 \item
\begin{equation}\label{metszesiclaim}
 F\in \overline{A_M} \hbox{ \ \ iff \ \ } \exists e\in M\; \exists f\in F \hbox{\ with \ } |e\cap f|=1.
\end{equation}
\item Matchings $M_1$ and $M_2$
 are in conflict iff $A_{M_1}\land A_{M_2}=\emptyset$.
\item 
If the matchings $F$ and $M$ are not in conflict, then 
\begin{equation}\label{metszesiclaim2}
\overline{A_{M\setminus F}}\subseteq  \overline{A_M} \hbox{\ \ and \ \ } \overline{A_M}\land A_F=\overline{A_{M\setminus F}}\land A_F.
\end{equation}
\end{enumerate}
\end{lemma}
\begin{theorem}\label{negdep}
  Let $\mathcal{M}$ be a collection of matchings in $K_N$ or $K_{N,N'}$.
The graph $G=G(\mathcal{M})$  described below is a negative dependency graph for the canonical
events $\{A_M\mid M\in \mathcal{M}\}$:
\begin{itemize}
\item $V(G)=\mathcal{M}$,
\item $E(G)=\Bigl\{\{M_1,M_2\}\mid M_1\in\mathcal{M}  \mbox{ and } M_2\in\mathcal{M}
 \mbox{ are in conflict}\Bigl\}.$
\end{itemize}
\end{theorem}

\noindent {\bf Proof:} For complete bipartite graphs we proved this theorem
in \cite{soda}, and therefore we have to prove it now for $K_N$. 
We will prove the theorem by induction on $N$. The base case 
$N=2$ is trivial.
Throughout this paper, we always assume that the vertex set of $K_N$ is
$[N]=\{1,2,\ldots, N\}$. There is a canonical injection from $[N]$ to $[N+s]$,
and consequently from 
$V(K_N)$ to $V(K_{N+s})$ and from $E(K_N)$ to $E(K_{N+s})$.
 Through this canonical injection,
every matching of $K_{N}$ can be viewed as a matching
of $K_{N+s}$. (Note that a perfect matching in $K_N$ will not be perfect
in $K_{N+s}$ for $s>0$.) {\sl To emphasize the difference in the size of the vertex set, we
use $A^N_M$ to denote the event  induced by the matching $M$ among the matchings 
of an $N$-vertex complete graph.}

\begin{lemma}\label{l1}
For any collection $\mathcal{M}$ of matchings in $K_N$,
 we have
\begin{equation}
  \nonumber
  \Pr(\wedge_{M\in\mathcal{M}}\overline{ A^N_M}) \leq
  \Pr(\wedge_{M\in\mathcal{M}}\overline{ A^{N+2}_M}).
\end{equation}
\end{lemma}

\noindent
{\bf Proof:} We partition the space of $\Omega_{N+2}$  into
$N+1$ sets as follows:
 for $1\leq i \leq N+1$, let $\mathcal{C}_i$ be the set of
perfect  matchings containing the edge $i(N+2)$.
We have
\begin{equation}
  \nonumber
  \Pr(\wedge_{M\in\mathcal{M}}\overline{ A^{N+2}_M})=\sum_{i=1}^{N+1}
\Pr(\wedge_{M\in\mathcal{M}}\overline{ A^{N+2}_M}\wedge \mathcal{C}_i).
\end{equation}
We observe that  $\mathcal{C}_i\subseteq \overline{ A^{N+2}_M}$ if and only if
$M$ conflicts $i(N+2)$, a one-edge matching. Let $\mathcal{B}_i$ be
 the subset of $\mathcal{M}$, whose elements are not in conflict with
the edge $i(N+2)$.
(In particular, $\mathcal{B}_{N+1}=\mathcal{M}$.)
We have
\begin{equation}
  \nonumber
  \wedge_{M\in\mathcal{M}}\overline{ A^{N+2}_M}\wedge \mathcal{C}_i
= \wedge_{M\in\mathcal{B}_i}\overline{ A^{N+2}_M}\wedge \mathcal{C}_i.
\end{equation}
Let $\phi_i$ be the transposition $i\leftrightarrow N+1$ acting on the set $\{1,2,...,N+2\}$.
Note that $\phi_i$ stabilizes $\mathcal{B}_{i}$,
interchanges $\mathcal{C}_i$ and $\mathcal{C}_{N+1}$, 
 and
maps $\wedge_{M\in\mathcal{B}_i}\overline{ A^{N+2}_M}\wedge \mathcal{C}_i$
to
$\wedge_{M\in\mathcal{B}_i}\overline{ A^{N+2}_M}\wedge C_{N+1}$.
We have
\begin{eqnarray}
   \Pr(\wedge_{M\in\mathcal{M}}\overline{ A^{N+2}_M})&=&\sum_{i=1}^{N+1}
\Pr(\wedge_{M\in\mathcal{M}}\overline{ A^{N+2}_M}\wedge \mathcal{C}_i)
 \label{kezdet} \\
&=&\sum_{i=1}^{N+1}
\Pr(\wedge_{M\in\mathcal{B}_i}\overline{ A^{N+2}_{M}}\wedge \mathcal{C}_i)
\nonumber  \\
&=&\sum_{i=1}^{N+1}
\Pr(\wedge_{M\in\mathcal{B}_i}\overline{ A^{N+2}_{M}}\wedge \mathcal{C}_{N+1})
\nonumber
 \\
&=& \sum_{i=1}^{N+1}
\Pr(\wedge_{M\in\mathcal{B}_i}\overline{ A^{N+2}_{M}} \mid \mathcal{C}_{N+1})
\Pr(\mathcal{C}_{N+1})\nonumber \\
&=&\frac{1}{N+1}\sum_{i=1}^{N+1}
\Pr(\wedge_{M\in\mathcal{B}_i}\overline{ A^{N}_{M}}),\label{veg}  \\
\hbox{\noindent \rm and estimating  further}\nonumber \\
&\geq & (N+1) \Pr(\wedge_{M\in\mathcal{M}}\overline{ A^{N}_M}) \frac{1}{N+1}
\nonumber\\
&=& \Pr(\wedge_{M\in\mathcal{M}}\overline{ A^{N}_M}). \nonumber
\end{eqnarray}
The proof of  Lemma \ref{l1} is finished. \hfill$\square$

For the completeness, we provide the variation of Lemma \ref{l1} for
the case of  $K_{N.N'}$. The proof will be omitted.
Let $A^{N,N'}_M$ be the event induced by the matching $M$ among the matchings 
of a complete bipartite graph $K_{N, N'}$.  
\begin{lemma}\label{l1'}
For any collection $\mathcal{M}$ of matchings in $K_{N,N'}$,
 we have
\begin{equation}
  \nonumber
  \Pr(\wedge_{M\in\mathcal{M}}\overline{ A^{N,N'}_M}) \leq
  \Pr(\wedge_{M\in\mathcal{M}}\overline{ A^{N+1, N'+1}}_M).
\end{equation}
\end{lemma}

\noindent {We are back to the  proof of Theorem \ref{negdep}:}
For any fixed matching $M\in \mathcal{M}$, and a subset
$\mathcal{J}\subseteq \mathcal{M}$ satisfying that for every
 $M'\in\mathcal{J}  $,
$M'$ is not in  conflict with $M$, by (\ref{eq:neg1}) it suffices to show that
\begin{equation} \label{cel}
  \Pr(\wedge_{M'\in \mathcal{J}} \overline{ A_{M'}}\mid A_M) \leq
\Pr(\wedge_{M'\in \mathcal{J}} \overline{ A_{M'}}).
\end{equation}
Let $\mathcal{J}'=\{M'\setminus M \mid M'\in \mathcal{J}\}$. 
Assume first that $\emptyset \notin \mathcal{J}'$.
Since every matching $M'$
in $\mathcal{J}$ is not in  conflict with $M$,
 the vertex set $V(M'\setminus M)$  of $M'\setminus M$ 
is disjoint from the vertex set $V(M)$ of $M$. 
Let $T=V(M)$ be the set of vertices covered by the matching $M$ and
$U$ be the set of vertices covered by at least one  matching 
$F\in \mathcal{J}'$.
We have $T\cap U=\emptyset$. Let $\pi$ be a permutation
of $[N]$ mapping $T$ to $\{N-|T|+1,N-|T|+2,\ldots, N\}$.
We have $\pi(U)\cap \pi(T)=\emptyset$. Thus, $\pi(U)\subseteq [N-|T|]$.
For a matching $F$, define another matching $\pi(F)$ by 
$\{\pi(u),\pi(v)\}\in \pi(F)$ if and only if  
$\{u,v\}\in F$. 
Let $\pi(\mathcal{J}')=\{\pi(F)\mid F\in\mathcal{J}' \}$ and $F'=\pi(F)$.
 Each matching in $\pi(\mathcal{J}')$ is
a matching in $K_{N-|T|}$. We obtain  (\ref{cel})
using Lemma \ref{l1} repeatedly:
\begin{eqnarray*}
   \Pr(\wedge_{M'\in \mathcal{J}} \overline{A_{M'}}\mid A_M)
&=& \frac{\Pr(\wedge_{M'\in\mathcal{J} } \overline{A_{M'}}\wedge A_M)}{\Pr(A_M)}\\
&=& \frac{\Pr(\wedge_{M'\in \mathcal{J}} 
\overline{A_{M'\setminus M}}\wedge A_M)}{\Pr(A_M)}  \hbox{\ by Lemma~\ref{onconflict}} \\
&=& \frac{\Pr(\wedge_{F\in \mathcal{J}'} \overline{A_F}\wedge A_M)}{\Pr(A_M)}\\
&=& \Pr(\wedge_{F\in \mathcal{J}'} \overline{A_F}\mid A_M)\\
&=&\Pr(\wedge_{F'\in \pi(\mathcal{J}')} \overline{A^N_{F'}}\mid A_{\pi(M)})\\
&=& \Pr(\wedge_{F'\in \pi(\mathcal{J}')} \overline{A^{N-|T|}_{F'}})\\
&\leq& \Pr(\wedge_{F'\in \pi(\mathcal{J}')} \overline{A^{N}_{F'}})   \hbox{\ by Lemma~\ref{l1}}   \\
&=& \Pr(\wedge_{F\in \mathcal{J}'} \overline{A^{N}_{F}})\\
&=& \Pr(\wedge_{M'\in \mathcal{J}} \overline{A^{N}_{M'\setminus M}})\\
&\leq& \Pr(\wedge_{M'\in \mathcal{J}} \overline{A^{N}_{M'}}).
\end{eqnarray*}
If  $\emptyset \in \mathcal{J}'$, then the LHS of the estimate above is zero,
and therefore we have nothing to do.
\hfill $\square$

The following example shows that in Theorem~\ref{negdep} one cannot have an
arbitrary graph in the place of $K_N$ or $K_{N,N'}$. Consider $G=C_6$, this graph
has two perfect matchings. Let $e$ and $f$ denote two opposite edges
of $C_6$. Consider the following two partial matchings: $\{e\}$ and
$\{f\}$. We have $\Pr(A_{\{e\}})= \Pr(A_{\{f\}})=  1/2$. However, we have
$$\Pr(A_{\{e\}}\mid \overline{A_{\{f\}}})=\frac{\Pr( A_{\{e\}}\wedge 
\overline{A_{\{f\}}})}{\Pr(\overline{A_{\{f\}}})} \not\leq \Pr( A_{\{e\}}).$$

Next, we prove a partial converse of Lemma~\ref{l1}.
\begin{lemma}\label{l6} Consider  a collection $\mathcal{M}$ of matchings in $K_N$, so that 
their canonical events satisfy condition (\ref{altcond}) for an $\epsilon < 1/4$, and in addition,
for any $uv\in E(K_N)$
\begin{equation}\label{parallel}
\sum_{M:uv\in M\in \mathcal{M}} \Pr(A^N_M)+ 2\Pr^2(A^N_M)<\epsilon.
\end{equation}
Then we have
\begin{equation}
  \nonumber
  \Pr(\wedge_{M\in\mathcal{M}}\overline{ A^{N+2}_M}) \leq
  y^2(\epsilon)
  \Pr(\wedge_{M\in\mathcal{M}}\overline{ A^{N}_M}).
\end{equation}
\end{lemma}
\noindent
{\bf Proof:} Partition $\Omega_{N+2} $, introduce ${\cal C}_i$ and 
${\cal B}_i$ as in the proof of Lemma \ref{l1}, and use the fact derived there
between (\ref{kezdet}) and (\ref{veg}) that 
\begin{equation} \label{sumcite}
   \Pr(\wedge_{M\in\mathcal{M}}\overline{ A^{N+2}_M})=
\frac{1}{N+1}\sum_{i=1}^{N+1}
\Pr(\wedge_{M\in\mathcal{B}_i}\overline{ A^{N}_{M}}).
\end{equation}
We are going to apply Theorem~\ref{simple} part (i) with  $S=\M\setminus \mathcal{B}_i$
and $T=\mathcal{B}_i$. $T=\mathcal{B}_i$ contains those matchings from $\mathcal{M}$ whose
support do not contain $i$, while $S$ contains those matchings whose support do contain $i$.
 We are going to show 
 \begin{equation} \label{applytoveg}
\frac{\Pr(\wedge_{M\in\M}\overline{ A^{N}_{M}})}{\Pr(
\wedge_{M\in\mathcal{B}_i}\overline{ A^{N}_{M}})}=
\Pr(\wedge_{M\in\M\setminus \mathcal{B}_i}\overline{ A^{N}_{M}} \mid
\wedge_{M\in\mathcal{B}_i}\overline{ A^{N}_{M}})\geq y(\epsilon)^{-2}.
\end{equation}
 We have from (\ref{l5formula})
\begin{equation} \label{applytoedges}
\Pr(\wedge_{M\in\M\setminus \mathcal{B}_i}\overline{ A^{N}_{M}} \mid
\wedge_{M\in\mathcal{B}_i}\overline{ A^{N}_{M}})\geq
\prod_{M\in   \mathcal{M}:i\in supp(M)  } \Bigl(1- \Pr(A_M)y(\epsilon)\Bigl).
\end{equation}
If the product in (\ref{applytoedges}) is empty, then 
we have nothing to prove in (\ref{applytoveg}). If there are $u\not= v$
such that $iu\in M_1\in \M$ and $iv\in M_2\in \M$, then $\{M\in   \mathcal{M}|  i\in supp(M)\}\subseteq
N_G(M_1)\cup N_G(M_2)$ (meaning neighborhoods in the conflict graph), and  the RHS of (\ref{applytoedges}) has a lower bound of 
$$
\prod_{M\in N_G(M_1)} \Bigl(1- \Pr(A_M)y(\epsilon)\Bigl)
\prod_{M\in N_G(M_2)} \Bigl(1- \Pr(A_M)y(\epsilon)\Bigl)\geq e^{-2\epsilon y(\epsilon)} =
y(\epsilon)^{-2},$$ like in the last line of the  proof of Theorem~\ref{simple}(ii), also using (\ref{altcond}).  
If there is an $ij$ edge, such that
$i\in supp(M)$ for $M   \in    \mathcal{M}$ implies $ij\in M$, then  condition 
(\ref{parallel}) gives a lower bound of $  y(\epsilon)^{-1}$ in a similar way for the RHS of
(\ref{applytoedges}).
We have from (\ref{sumcite}) and the estimate above:
\begin{eqnarray*}
\Pr(\wedge_{M\in\mathcal{M}}\overline{ A^{N+2}_M})
&=& \frac{1}{N+1}\sum_{i=1}^{N+1}
\Pr(\wedge_{M\in\mathcal{B}_i}\overline{ A^{N}_{M}}) \\
&\leq& \frac{1}{N+1}\sum_{i=1}^{N+1}
\Pr(\wedge_{M\in\mathcal{M}}\overline{ A^{N}_{M}})
  y^2(\epsilon)\\
&=&  y^2(\epsilon)
\Pr(\wedge_{M\in\mathcal{M}}\overline{ A^{N}_M}).
\end{eqnarray*}
The proof of Lemma  \ref{l6} is finished. \hfill$\square$

Here is a similar Lemma  for $K_{N,N'}$. The proof is similar and
will be omitted. 
\begin{lemma}\label{l6'} Consider  a collection $\mathcal{M}$ of matchings in $K_{N,N'}$, so that 
their canonical events satisfy condition (\ref{altcond}) for an $\epsilon < 1/4$, and in addition,
for any $uv\in E(K_{N,N'})$
\begin{equation}\label{parallel'}
\sum_{M:uv\in M\in \mathcal{M}} \Pr(A^{N,N'}_M)+ 2{\Pr}^2(A^{N,N'}_M)<\epsilon.
\end{equation}
Then we have
\begin{equation}
  \nonumber
  \Pr(\wedge_{M\in\mathcal{M}}\overline{ A^{N+1,N'+1}_M}) \leq
  y^2(\epsilon)
  \Pr(\wedge_{M\in\mathcal{M}}\overline{ A^{N,N'}_M}).
\end{equation}
\end{lemma}

\section{Upper bounds in the matching models}

Now we consider $\Omega$, the uniform probability space of perfect matchings
in $K_N$ ($N$ even) or $K_{N,N'}$ (with $N\leq N'$). 
 Let $\mathcal{M}$ be a collection of partial matchings. 
For any $F\in \M$, 
let 
$$\M_F=\{M\setminus F\mid M\in \M, M\not=F, M\cap F\not=\emptyset, F
\mbox{ is not in conflict to } M\}.$$ 

We say that a collection of matchings $\M$ in $K_N$
is {\em $\delta$-sparse} if
\begin{enumerate}
\item No matching from $\M$ is a subset of another matching from $\M$.
\item 
  $\M$ satisfies (\ref{altcond}) and (\ref{parallel}) with $\delta$ instead of $\epsilon$.
\item For any $F\in \M$, 
\begin{equation}\label{eq:23}
\sum_{H:H\in \M_F} \Pr(A_H^{N-2|F|})+\Pr(A_H^{N-2|F|})^2 <\delta,
\end{equation}
where $A_H^{N-2|F|}$ indicates that vertices of $F$ has been removed from
the underlying vertex set $[N]$ when creating $\Omega$.
\end{enumerate}
Similarly, a collection of matchings $\M$ in $K_{N,N'}$ is {\em
  $\delta$-sparse} if  
\begin{enumerate}
\item No matching from $\M$ is a subset of another matching from $\M$.
\item 
  $\M$ satisfies (\ref{altcond}) and (\ref{parallel'}) with $\delta$ instead of $\epsilon$.
\item For any $F\in \M$, 
\begin{equation}\label{eq:23'}
\sum_{H:H\in \M_F} \Pr(A_H^{N-|F|, N'-|F|})+\Pr(A_H^{N-|F|,N'-|F|})^2 <\delta,
\end{equation}
where $A_H^{N-|F|, N'-|F|}$ indicates that vertices of $F$ has been removed from
the vertex set of $K_{N,N'}$ when creating $\Omega$.
\end{enumerate}

For a positive integer $r$, we say that $\M$ is {\em $r$-bounded},  if for all $M\in\M$, 
$|M|\leq r$.

The main result of this section is  the following theorem.
\begin{theorem}\label{posdep}
  Let $\mathcal{M}$ be a collection of matchings in $K_N$ or $K_{N,N'}$.
If $\M$ is $\delta$-sparse and $r$-bounded, then the negative dependency graph is also
an $\epsilon$-near-positive dependency graph with
\begin{equation} \label{elsoresz}
\epsilon=1-e^{-\delta y(2\delta)-
\delta^2y^2(2\delta)}y^{-2r}(2\delta)
\end{equation}
and therefore
\begin{equation} \label{masresz}
\Pr(\wedge_{M\in \M} \overline{ A_M})
\leq \prod_{M\in \M}\Biggl(1-\Pr(A_M)e^{-\delta y(2\delta)-
\delta^2y^2(2\delta)}y^{-2r}(2\delta)
\Biggl).
\end{equation}
\end{theorem}
We are going to prove  Theorem~\ref{posdep} for $K_N$, and leave  the proof
for  $K_{N,N'}$, which requires only negligible changes,  to the
Reader. (More explicitly, one need to replace $A^{N}_M$ to $A^{N,N'}_M$,
$A^{N+2}_M$ to $A^{N+1,N'+1}$, $A^{N-2|F|}_M$ to $A^{N-|F|,N'-|F|}_M$, 
and Lemma \ref{l6} to Lemma \ref{l6'}.)

\noindent
{\bf Proof of Theorem \ref{posdep}:}
We are going to
 show that the negative dependency graph $G$ defined for matchings 
of $K_N$ in $\M$  is also
an $\epsilon$-near-positive dependency graph with $\epsilon$ as in 
  (\ref{elsoresz}); 
and then Theorem ~\ref{ub} together with  (\ref{elsoresz})
 will finish the proof of (\ref{masresz}) and Theorem ~\ref{posdep}.
The first part of the definition, $\Pr(A_i\wedge A_j)=0$ for $ij$ edges 
comes for free. We focus on the second part.

For any $F\in \M$ and a subset $S\subseteq \overline{ N_G(F)}$,
we need to prove
\[\Pr(A_{F}\mid \wedge_{M\in S} \overline{ A_M})\geq (1-\epsilon) \Pr(A_{F}),\]
or equivalently,
\[\Pr(\wedge_{M\in S} \overline{ A_M} \mid A_F) 
\geq (1-\epsilon)\Pr(\wedge_{M\in S} \overline{ A_M}).\]
Let $S_F=\{M\setminus F\mid M\in S\}$. Observe that $\emptyset \notin S_F$. 
Note that
\begin{eqnarray}
\Pr(\wedge_{M\in S} \overline{ A_M} \mid A_F)
&=&\frac{\Pr(\wedge_{M\in S} \overline{ A_M}\wedge  A_F)}{\Pr(A_F)}
\label{szam1}\\
&=& \frac{\Pr(\wedge_{M\in S} \overline{ A_{M\setminus F}}\wedge  A_F)}
{\Pr(A_F)} \nonumber \\
&=&  \Pr(\wedge_{M\in S_F} \overline{ A_M} \mid A_F). \label{szam1veg}
\end{eqnarray}
We have
\begin{eqnarray}
 \Pr(\wedge_{M\in S_F} \overline{ A_M} \mid A_F)
&=& \Pr(\wedge_{M\in S_F} \overline{ A^{N-2|F|}_M}) \label{szam3}\\
&=& \Pr(\wedge_{M\in S_F} \overline{ A^{N}_M})\prod_{j=1}^{|F|} \nonumber
\frac{\Pr(\wedge_{M\in S_F} \overline{ A^{N-2j}_M})}
{\Pr(\wedge_{M\in S_F} \overline{ A^{N-2j+2}_M})}\\
\hbox{\rm (by Lemma \ref{l6})} 
&\geq&  \Pr(\wedge_{M\in S_F} \overline{ A^{N}_M}) \nonumber
\prod_{\ell=0}^{|F|-1} y^{-2}(2\delta)\\
&\geq&  \Pr(\wedge_{M\in S_F} \overline{ A^{N}_M}) 
y^{-2r}(2\delta).
\label{szam3end}
\end{eqnarray}
(Note that condition (\ref{parallel}) is implied by assumption 3.) 
For any $M$, which does not conflict to $F$, 
we have $\overline{ A_{M\setminus F}}\subset \overline{ A_M}$.
We have with $S_F=\{M\setminus F\mid M\in S\}$ that
\begin{eqnarray}
 \frac{\Pr(\wedge_{M\in S_F}\overline{ A^{N}_M})}{\Pr(\wedge_{M\in S}
\overline {A^{N}_M})} \label{szam2}
&=& \frac{\Pr(\wedge_{M\in S}\overline{ A^{N}_{M\setminus F}})}
{ \Pr(\wedge_{M\in S}\overline{ A^{N}_M})}\\
&=& \frac{\Pr(\wedge_{M\in S}\overline{ A^{N}_{M\setminus F}}\wedge
 \overline{ A^{N}_M}) \nonumber
}{ \Pr(\wedge_{M\in S}\overline{ A^{N}_M})}\\
&=& \frac{\Pr([\wedge_{{M\in S}, {M\cap F\not=\emptyset}}
\overline{ A^{N}_{M\setminus F}}]\wedge \nonumber
[\wedge_{M\in S} \overline{ A^{N}_M}])}{ \Pr(\wedge_{M\in S}\overline{ A^{N}_M})}\\
&=& \Pr(\wedge_{M\in S_F\setminus S}\overline{ A^{N}_M}\mid \wedge_{M\in S}
\overline{ A^{N}_M}). \label{szam2end}
\end{eqnarray}
Now apply Theorem~\ref{simple} part (i) to $S_F\setminus S$, $S$ 
and 
 $S\cup S_F$ instead of  $S$, $\T$ and   $\M$:
 \begin{eqnarray}\nonumber
& & \Pr\Bigl(\land_{M\in S_F\setminus S}\overline{A_M^N}\mid \land_{M\in S}\overline{A_M^N}\Bigl)
 \geq \prod_{M\in S_F\setminus S} \bigl(1-\Pr({A_M^N})y(2\delta)\bigl)\\
& \geq & \exp\Bigl(-\sum_{M\in S_F\setminus S} \Pr({A_M^N})y(2\delta) \nonumber
-\sum_{M\in S_F\setminus S} \Pr({A_M^N})^2y^2(2\delta)\Bigl)\\
&\geq & e^{-\delta y(2\delta)- \label{delta2}
\delta^2y^2(2\delta)}. 
 \end{eqnarray}
 Finally, we have
\begin{eqnarray*}
 & & \Pr(\wedge _{M\in S} \overline{ A_M}\mid A_F)\\
\hbox{\rm by (\ref{szam1}-\ref{szam1veg})} &=& 
\Pr(\wedge _{M\in S_F} \overline{ A_M}\mid A_F) \\
\hbox{\rm by (\ref{szam3}-\ref{szam3end})}&\geq& 
\Pr(\wedge_{M\in S_F} \overline{ A^{N}_M}) y^{-2r}(2\delta)\\
\hbox{\rm by (\ref{szam2}-\ref{szam2end})}
&=& \Pr(\wedge_{M\in S} \overline{ A^{N}_M}) 
\Pr(\wedge_{M\in S_F\setminus S}\overline{ A^{N}_M}\mid \wedge_{M\in S}
\overline{ A^{N}_M})y^{-2r}(2\delta)  \\
\hbox{\rm by (\ref{delta2}) }&\geq& \Pr(\wedge_{M\in S} \overline{ A^{N}_M}) 
e^{-\delta y(2\delta)-
\delta^2y^2(2\delta)}y^{-2r}(2\delta).
\end{eqnarray*}
Thus, the negative dependency graph $G$ is also a $\epsilon$-near-positive
dependency graph. The proof is finished by Theorem \ref{ub}.
\hfill $\square$

Theorem \ref{simple} provides a lower bound on $\Pr(\wedge_{M\in \M}\overline{A_M})$
while Theorem \ref{posdep} provides an upper bound on $\Pr(\wedge_{M\in \M}\overline{A_M})$.
Under  proper conditions, the combination of the two theorems gives asymptotics
for $\Pr(\wedge_{M\in \M}\overline{A_M})$, like in the following theorem.

\begin{theorem} \label{enum} Let $\Omega$ be the uniform probability space of perfect matchings
in $K_N$ ($N$ even) or $K_{N,N'}$ (with $N\leq N'$). 
Let   $r=r(N)$ be a positive integer and 
$1/16>\epsilon=\epsilon(N)>0$
as $N$ approaches  infinity.
Let $\M=\M(N)$ be a collection  of matchings in $K_N$ or $K_{N,N'}$, respectively,
such that none of these matchings is a subset of another.
For any $M\in \M$, let $A_M$ be the event consisting of perfect matchings 
extending $M$.
 Set $\mu=\mu(N)=\sum_{M\in\M}\Pr(A_M)$. Suppose that $\M$ satisfies
\begin{enumerate}
\item $|M|\leq r$, for each $M\in \M$.
\item $\Pr(A_M)<\epsilon$ for each $M\in \M$.
\item $\sum_{M':A_{M'}\cap A_M=\emptyset}\Pr(A_{M'})<\epsilon$ for each $M\in \M$.
\item $\sum_{M: uv\in M\in\M}\Pr(A_{M})<\epsilon$ for each single edge $uv$.
\item $\sum_{H\in \M_F}\Pr(A_H^{N-2r})<\epsilon$ for each $F\in \M$.
\end{enumerate}
Then, if $r\epsilon=o(1)$, we have
\begin{equation}
  \label{eq:enum}
\Pr(\wedge_{M\in \M}\overline{A_M})= e^{-\mu +O(r\epsilon \mu)}, 
\end{equation}
and furthermore, if $r\epsilon \mu=o(1)$, then 
\begin{equation}
  \label{eq:enum2}
\Pr(\wedge_{M\in \M}\overline{A_M})= \Bigl(1+ O(r\epsilon \mu)\Bigl)e^{-\mu}.
\end{equation}
\end{theorem}
{\bf Proof:} Let $G$ be the graph defined in Theorem \ref{negdep}.
By Theorem \ref{negdep}, the graph $G$ is a negative dependency graph
for the family of canonical events $\{A_M\}_{M\in \M}$. 
Note that the condition (\ref{altcond}) in Theorem \ref{simple}
is satisfied with $2\epsilon$, where $\epsilon$ is from the conditions of Theorem~\ref{enum}, instead of $\epsilon$. Applying Theorems \ref{simple} and  \ref{negdep},
we have 
\begin{eqnarray*}
  \Pr(\wedge_{M\in \M}\overline{A_M}) &\geq& 
\exp\Biggl(-\sum_{M\in \M}\Pr(A_M) y(2\epsilon)-   \sum_{M\in \M}\Pr^2(A_M) y^2(2\epsilon)\Biggl)\\
&>& \exp \Biggl(-\sum_{M\in \M}\Pr(A_M) y(2\epsilon)-\sum_{M\in \M}\Pr(A_M)\epsilon y^2(2\epsilon)\Biggl)\\
&=& \exp \Bigl(-\mu\bigl( 1+ 3\epsilon +O(\epsilon^2)\bigl) \Bigl).
\end{eqnarray*}

Now we consider the upper bound. Note that $\M$ is $2\epsilon$-sparse and $r$-bounded.
By Theorem \ref{posdep},
we have 
\begin{eqnarray*}
   \Pr(\wedge_{M\in \M}\overline{A_M}) &\leq& \prod_{M\in \M}\Bigl(1-\Pr(A_M)e^{-2\epsilon y(4\epsilon)-
(2\epsilon)^2y^2(4\epsilon)}y^{-2r}(4\epsilon) \Bigl)\\
&\leq& \exp\Biggl(- \sum_{M\in \M} \Pr(A_M)e^{-2\epsilon y(4\epsilon)-
(2\epsilon)^2y^2(4\epsilon)}y^{-2r}(4\epsilon) \Biggl)\\
&=& \exp\Bigl(-\mu \bigl(1-(8r+2)\epsilon + O(r^2\epsilon^2)\bigl)\Bigl).
\end{eqnarray*}
Combining the lower bound and the upper bound above, we obtain
equation (\ref{eq:enum}). 
\hfill $\square$

\section{Asymptotic results in the matching models}

\subsection{Applications I: Counting $k$-cycle free permutations and
Latin rectangles}
It is well-known and easy to see
that for any fixed $k$, the probability that a random permutation $\pi\in S_N$
is $k$-cycle free is   $\sim e^{-1/k}$, see 
 \cite{wilf} or \cite{borkowitz}. Earlier   we \cite{soda} obtained an
$(1-o(1))e^{-1/k}$ lower bound for this probability from the
Lov\'asz Local Lemma. To illustrate the applicability of Theorem~\ref{enum},
we show a lesser known result: the very same asymptotic formula  holds whenever $k=o(N)$. 

 Let us be given two $N$-element sets with elements $\{1,2,...,N\} $ and \linebreak
$\{1',2',...,N'\} $. Let us identify a permutation of the first set,
$\pi$, with a matching between the two sets, such that $i$ is joined to
$\pi(i)'$. A $k$-cycle in the permutation can be identified with a
matching between $K\subset \{1,2,...,N\}$    (with $|K|=k$)   and  $\{\ell':\ \ell\in K\}$,
which does not have a proper non-empty  subset $K_1\subset K$, such that  
the matching also matches $K_1$ to $\{\ell':\ \ell\in K_1\}$.
The bad events for the negative dependency graph
are these $k$-element matchings; there are ${N\choose k}(k-1)!$ of them.
We have $|\M|={N\choose k}(k-1)!$.  For each $M\in \M$, we have
$\Pr(A_M)=\frac{1}{{N\choose k} k!}$. Two matchings, $M,M'
\in \M$, $M\not= M'$, conflict  each other if and only if the two cycles have non-empty intersection, i.e. have common elements.

Let $r=k$ and $\epsilon=\frac{k}{N-k+1}$. Now we will verify the conditions of 
Theorem \ref{enum}. Items 1 and 2 are satisfied by our choice of $r$ and $\epsilon$.
For item 3, we have
\begin{eqnarray}
  \sum_{M':A_{M'}\cap A_M=\emptyset}\Pr(A_{M'}) &=& 
\left( {N\choose k}(k-1)! - {N-k\choose k}(k-1)! \right) \frac{1}{{N\choose k} k!}\nonumber\\
&=& \frac{1}{k} \Biggl(1- \prod_{i=1}^k \frac{N-k-i+1}{N-i+1}\Biggl)\nonumber\\
&=& \frac{1}{k} \Biggl(1- \prod_{i=1}^k \left(1- \frac{k}{N-i+1}\right)\Biggl)\nonumber\\
&<& \frac{1}{k} \sum_{i=1}^k \frac{k}{N-i+1}\leq \frac{k}{N-k+1}= \epsilon.\label{lastlabel}
\end{eqnarray}
Now we verify item 4. For any $uv\in M\in \M$, a $k$-matching $M$
contains a given edge $uv$, if and only if $v=\pi(u)'$ for some
$k$-cycle permutation $\pi$.  The number of such $k$-cycles is
${N\choose k-2}(k-2)!$. We have
\begin{eqnarray*}
  \sum_{M: uv\in M\in\M}\Pr(A_{M})&=& {N\choose k-2}(k-2)!\frac{1}{{N\choose k} k!}\\
&=& \frac{1}{(N-k+2)(N-k+1)}<\epsilon.
\end{eqnarray*}
For any $F\in \M$, now $\M_F$ is empty in our
special setting, hence item 5 holds trivially. All conditions of Theorem \ref{enum} 
are verified. Observe 
\begin{equation}\label{muu}
\mu=\sum_{M\in\M}\Pr(A_M)= {N\choose k}(k-1)!\frac{1}{{N\choose k} k!}=\frac{1}{k}.
\end{equation}
 Therefore Theorem~\ref{enum}
applies, and the number of $k$-cycle-free permutations is $(1+O(k/N))e^{-1/k}$.
\cite{wilf} goes further than this, and gives asymptotic formula for the number of permutations without cycles of length $r$ or less, for fixed $r$. Simple generating 
function arguments would not allow $k$ (or $r$) to be variables. However, our method
allows it. The following result first occured in \cite{bender}:
\begin{theorem}
Let us be given a $K\subset\{1,2,...,N\}$ and set $r=\max K$. Assume that 
$$r^2\Bigl(\sum_{k\in K} \frac{1}{N-k+1}
\Bigl)\rightarrow 0,  \hbox{\ \ \ and \ \ \ }  R=r^2\Bigl(\sum_{k\in K} \frac{1}{k}\Bigl)\Bigl(\sum_{k\in K} \frac{1}{N-k+1}
\Bigl)\rightarrow 0.$$
Then, the probability that a random permutation of $N$ elements do not contain
any cycle, whose length belongs to $K$, is $\displaystyle \bigl(1+O(R)\bigl)\exp\Bigl(-\sum_{k\in K} \frac{1}{k}\Bigl)$.
\end{theorem}
\noindent {\bf Proof:} The proof above goes through with minor modifications. Set $\epsilon=r\sum_{k\in K}
\frac{1}{N-k+1}$, change (\ref{muu}) to $\mu=
\sum_{k\in K}{N\choose k}(k-1)!\frac{1}{{N\choose k} k!}=\sum_{k\in K}\frac{1}{k}$,
and for a matching $M$ corresponding to an $\ell$-cycle,  change (\ref{lastlabel}) 
for the estimation of  $\sum_{M':A_{M'}\cap A_M=\emptyset}\Pr(A_{M'})$   to
$\sum_{k\in K}\left( {N\choose k}(k-1)! - {N-k\choose k}(k-1)! \right)
 \frac{1}{{N\choose k} k!}$ $
=\sum_{k\in K} \frac{1}{k} \Bigl(1- \prod_{i=1}^k \frac{N-\ell-i+1}{N-i+1}\Bigl)$
$=\sum_{k\in K}\frac{1}{k}
 \Bigl(1- \prod_{i=1}^k \left(1- \frac{\ell}{N-i+1}\right)\Bigl)$\\
  $<\sum_{k\in K} \frac{1}{k} \sum_{i=1}^k \frac{\ell}{N-i+1}\leq \sum_{k\in K}\frac{\ell}{N-k+1}\leq \epsilon.$
 \hfill$\square$
 \smallskip

Let us turn now to the enumeration of Latin rectangles. A $k\times n$
Latin rectangle is a sequence of $k$ permutations of $\{1,2,...,n\} $
written in a matrix form, such that no column has any repeated entries.
Let $L(k,n)$ denote the number of $k\times n$ Latin rectangles.
$L(2,n)$ is just $n!$ times the number of derangements, i.e. $(n!)^2e^{-1}$.
 In 1944, Riordan
\cite{riordan} showed that $L(3,n)\sim (n!)^3e^{-3}$. In 1946, 
Erd\H os and Kaplansky \cite{erdoskap} showed 
\begin{equation} \label{yam}
L(k,n)\sim (n!)^ke^{-{k\choose 2}}
\end{equation}
for $k=o\bigl((\log n)^{3/2}\bigl)$. In 1951, Yamamoto \cite{yamamoto} 
 extended this asymptotic 
formula for $k=o(n^{1/3})$.
In 1978, Stein \cite{stein} refined the asymptotic formula to 
\begin{equation} \label{latinas}
L(k,n)\sim (n!)^ke^{-{k\choose 2}-{k^3\over 6n}}
\end{equation} 
using  the Chen-Stein
method \cite{chen}, and extended the range to $k=o(n^{1/2})$.
The current best asymptotic formula is due to Godsil and McKay \cite{godsil},
whose further refined formula,  
 $L(k,n)\sim (n!)^k\Bigl(\frac{(n)_k}{n^k}\Bigl)^n\Bigl(1- \frac{k}{n}    \Bigl)^{-n/2}e^{-k/2}  $ works for $k=o(n^{6/7})$.

Formula (\ref{latinas}) has had an unexpected proof by Skau \cite{skau},
who proved, for any $1\leq k\leq n$, the inequality 
\begin{equation} \label{lat2}
(n!)^k\prod_{t=1}^{k-1} \Biggl(1-\frac{t}{n} \Biggl)^n\leq L(k,n)
\end{equation} 
from the van der Waerden inequality for the permanent. 
If $k=o((n/\log n)^{1/2})$, the lower bound in (\ref{lat2})  is asymptotically the same 
as the RHS of (\ref{latinas}). Skau's asymptotically tight upper bound \cite{skau}   followed
from Minc's inequality for the permanent.

In \cite{soda} we claimed   (\ref{lat2}) from the Lov\'asz Local Lemma in error.
However, our method still gives back Yamamoto's range for (\ref{yam}). 
 Fix an arbitrary $t\times n$ Latin rectangle  with rows
$\pi_1,\pi_2,...,\pi_t$. Consider a complete bipartite graph with classes
$\{1,2,...,n\} $ an $\{1',2',...,n'\} $, and let $\Omega $ be the space of perfect 
matchings in this  complete bipartite graph. Permutation $ \pi_{t+1}$ of $\{1,2,...,n\} $
are in one-to-one correspondence with perfect matchings by 
$(\pi_{t+1}(j),j'):1\leq j \leq n$. Permutation $ \pi_{t+1}$ fails to extend the given 
Latin rectangle into a $(t+1)\times n$ Latin rectangle iff there are $i,j$ such that
$\pi_i(j)=\pi_{t+1}(j)$. Therefore a perfect matching provides a legal $(t+1)^{th}$
row for the Latin rectangle iff it does not contain any of the edges 
$(\pi_{i}(j),j'):1\leq j \leq n, 1\leq i\leq t$.
Define the event $A_{ij} $ as the canonical event in $\Omega$ 
corresponding to the one-edge
matching $(\pi_{i}(j),j')$. 
Let $G$ be the   a negative dependency graph for the family of events
 $A_{ij}$, according to
Theorem \ref{negdep}. $G$ is $(t-1)$-regular.
We can apply Theorem~\ref{enum} with $\epsilon=2t/n$  and $\mu=\frac{1}{n}\cdot (nt)=t$. 
Condition  1 of Theorem~\ref{enum} holds with $r=1$, Condition 2 holds as $1/n<\epsilon$,
Condition  3 holds as   
$2(t-1)/n<\epsilon$, Condition 4 holds  like  Condition 2, and Condition 5 holds vacuously.
Hence
$\#\pi_{t+1}/n!=\exp\Biggl({-t+O\Bigl(\frac{t^2}{n}\Bigl)}\Biggl)$ by (\ref{eq:enum}),
and $L(k,n)=\prod_{t=0}^{k-1}n!\exp\Biggl({-t+O\Bigl(\frac{t^2}{n}\Bigl)}\Biggl)=
(n!)^k\exp\Biggl(-{k\choose 2}+O\Bigl(\frac{k^3}{n}\Bigl)\Biggl).$

\subsection{Applications II: The configuration model and the enumeration of
 $d$-regular graphs} \label{dreg}
For a given   sequence of positive integers with an even sum,
${\mathbf d} = (d_1, d_2, \ldots, d_n)$,
the {\em configuration model of random multigraphs
with degree sequence $\mathbf d$} is defined as follows \cite{configmodel}.

\noindent
1.\hspace{2mm} Let us be given a set $U$ that contains $N=\sum_{i=1}^n d_i $ 
distinct mini-vertices. Let $U$ be partitioned into $n$ classes such
 that the  $i$th
class consists of $d_i$ mini-vertices. This $i$th class   will be associated
with vertex $v_i$ after identifying its elements through a {\it projection}.\\
2.\hspace{2mm} Choose a random perfect matching $M$ of the mini-vertices in $U$ 
uniformly.\\
3.\hspace{2mm}  Define a random multigraph $G$ associated with $M$ as
follows: For any two (not necessarily distinct) 
vertices $v_i$ and $v_j$, the number of edges
joining $v_i$ and $v_j$ in $G$ is equal to the total number of edges
in $M$ between mini-vertices associated with $v_i$ and mini-vertices
associated with $v_j$. 

The configuration model of random $d$-regular graphs
on $n$ vertices is the instance $d_1= d_2= \cdots = d_n=d$, where $nd$ is even.

The enumeration problem of labelled $d$-regular graphs
 has a rich history in the literature. The first result was Bender and Canfield
\cite{benderca}, who showed in 1978 that for any fixed $d$, with $nd$  even,
the number of them is 
\begin{equation}\nonumber
\bigl(\sqrt{2}+o(1)\bigl)e^{\frac{1-d^2}{4}}\Biggl(\frac{d^dn^d}{e^d(d!)^2}\Biggl)^{\frac{n}{2}}.
\end{equation}
The same result was discovered at the same time by Wormald.
In 1980, Bollob\'as  \cite{configmodel}
introduced probability to this enumeration problem by
defining the configuration model, and put the result in the alternative form
\begin{equation} \label{os}
(1+o(1)) e^{\frac{1-d^2}{4}}\frac{(dn-1)!!}{(d!)^n}.
\end{equation}
where the term
 $(1+o(1))e^\frac{1-d^2}{4}$
in (\ref{os}) can be explained as the probability of obtaining a simple graph
after the projection. The semifactorial $(dn-1)!!=\frac{(dn)!}{(dn/2)! 
2^{dn/2}}$ equals  the number of perfect matchings on $dn$ 
elements, and $\frac{1}{(d!)^n}$ is just the number of ways matchings
can yield the same simple graph after projection. Bollob\'as also extended
the range of the asymptotic formula to $d<\sqrt{2\log n}$, which was
further extended to $d=o(n^{1/3})$ by  McKay \cite{72} in 1985.
The strongest result is due to McKay and Wormald \cite{81} in 1991,
who refined the probability of obtaining a simple graph
after the projection to 
\begin{equation}\label{strongest}
(1+o(1))e^{\frac{1-d^2}{4}-\frac{d^3}{12n}+O(\frac{d^2}{n})}
\end{equation}
 and extended
the range of the asymptotic formula to $d=o(n^{1/2})$. Wormald's Theorem 2.12
in \cite{wormald} (originally published in  \cite{wormaldos})
  asserts that for any fixed numbers $d\geq 3$
 and $g\geq 3$, the number
of labelled $d$-regular graphs with girth at least $g$, is
\begin{equation}\label{faktor}
(1+o(1))e^{-\sum_{i=1}^{g-1}\frac{(d-1)^i}{2i}}\frac{(dn-1)!!}{(d!)^n}.
\end{equation}
In our theorem below, we allow both $d$ and $g$ go to infinity slowly.
If we set $g=3$, we get back the asymptotic formula for the number of
$d$-regular graphs up to $d=o(n^{1/3})$, giving an alternative proof to McKay's
result cited above. However, our method inherently fails to extend this
result as McKay and Wormald did. In fact, our method already
fails to extend the lower bound.  McKay, Wormald and Wysocka
\cite{wysocka}  proved  Theorem~\ref{multigraph} below under a slightly weaker assumption
$(d-1)^{2g-3}=o(n)$. We could somewhat reduce the exponent in $g^6$ \cite{JCMCC}, but at least
a factor of $g$ comes from comes from the condition $r\epsilon\mu\rightarrow 0$ among
the conditions of Theorem~\ref{enum} that we use.  A power of $g$ in (\ref{condition}) 
is of secondary importance beside the exponential term.

\begin{theorem} \label{multigraph}
  In the configuration model, assume $d\geq 3$ and 
\begin{equation} \label{condition}
g^6(d-1)^{2g-3}=o(n).
\end{equation}
Then
the probability that the random $d$-regular multigraph has
girth at least $g\geq 1$ is $  (1+o(1))\exp\Bigl({-\sum_{i=1}^{g-1}\frac{(d-1)^i}{2i}} \Bigl) $,
and hence the number of $d$-regular graphs on $n$ vertices with girth
 at least $g\geq 3$ is
\[(1+o(1))e^{-\sum_{i=1}^{g-1}\frac{(d-1)^i}{2i}}\frac{(dn-1)!!}{(d!)^n}.\]
(The case $g=3$ means that the random $d$-regular multigraph is actually
a simple graph.) Furthermore, the number of $d$-regular graphs not containing 
cycles whose length is in a set ${\cal C}\subseteq\{3,4,...,g-1\}$, is
 \[(1+o(1))e^{- \frac{d-1}{2}-\frac{(d-1)^2}{4}  -\sum_{i\in {\cal C}}\frac{(d-1)^i}{2i}
 }\frac{(dn-1)!!}{(d!)^n}
. \]
\end{theorem}
{\bf Proof:} We prove the first claim. To prove the second claim, only
(\ref{errorinmainterm}) has to be adjusted, everything else remains the same.
For $i=1,2,\ldots, g-1$, let $\M_i$ be the set of  partial matchings of $U$ whose
projection gives  precisely a cycle of length $i$; there are {\em exactly}
$\frac{1}{2i}{n\choose i} i! d^i (d-1)^i$
of them.  The bad events for the negative dependency
graph are the union of matchings $\M=\cup_{i=1}^{g-1} \M_i$.
For each $M_i\in \M_i$ ($i=1,2,\ldots, g-1$), we have
\begin{equation}
  \label{eq:mi}
  \Pr(A_{M_i})=\frac{1}{(nd-1)(nd-3)\cdots (nd-2i+1)}.
\end{equation}
We have
\begin{eqnarray}
\sum_{M\in \M}\Pr(A_M) &=& \sum_{i=1}^{g-1} \frac{1}{2i}{n\choose i} i! d^i (d-1)^i
\nonumber
\frac{1}{(nd-1)(nd-3)\cdots (nd-2i+1)}\\
&=&  \sum_{i=1}^{g-1} \frac{(d-1)^i}{2i}\left(1+ O\left(\frac{i^2}{n}\right)\right)
= \left(1+ O\left(\frac{g^2}{n}\right)\right)   \sum_{i=1}^{g-1} \frac{(d-1)^i}{2i}   
\nonumber .\\ \label{errorinmainterm}
\end{eqnarray}

Let $r=g-1$ and $\epsilon=\frac{K'g^5(d-1)^{g-2}}{n}$ for a large constant  $K'$.
Now we verify the conditions of Theorem \ref{enum}. Item 1 and 2 are trivial
by the definition of $r$ and $\epsilon$.
Item 3 can be verified as follows. Consider now a fixed   $M\in \M_1$, then $M=\{e\}$, where $e=xy$, where $x$ and $y$ belong to the same
$d$-element $C$ that projects to a mini-vertex. If  $M'\in \M_1$ and $A_{M'}\cap A_M=\emptyset$, then $M'=\{f\}$, both endpoints of $f$ are in $C$,
and one of them is $x$ or $y$. We have exactly $2(d-2) $ of such $M'$ matchings, and for each $\Pr(A_{M'})=\frac{1}{nd-1}$. 
If $M'\in \M_i$ for some $i\geq 2$, then we see $i$ classes projecting to distinct mini-vertices, one of them is $C$, the remaining $i-1$ are arbitrary among the remaining
$n-1$. There are $i$ vertex disjoint edges between these classes, so that so that after the projection we see an $i$-cycle, and one of those $i$ edges has either $x$ or $y$ as an
endpoint. To build all such $M'$ matchings, select the $i-1$ classes in $\binom{n-1}{i-1}$ ways, put them into a cycle in $(i-1)!/2$ ways, select whether $x$ or $y$ will be an
endpoint of one of the $i$ edges in 2 ways, then select $x$'s neighbor in the class dictated by the cycle in $d$ ways, the endpoint of the next edge in $d-1$ ways, and keep going.   
When we return to $C$, we have  $d-1$ choices, as we cannot return to  the vertex of $\{x,y\}$, from which we started.   In addition, we have $\Pr(A_{M'})=
\frac{1}{(nd-1)(nd-3)\cdots (nd-2i+1)}$.
In conclusion, we obtain
\begin{eqnarray}
    \sum_{M':A_{M'}\cap A_M=\emptyset}\Pr(A_{M'}) &\leq& \nonumber
\frac{2d-4}{nd-1}+ 
\sum_{i=2}^{g-1} \sum_{M'\in \M_i: A_{M'}\cap A_M=\emptyset}\Pr(A_{M'})\\
\nonumber
&\leq& \frac{2d-4}{nd-1}+ \sum_{i=2}^{g-1} \frac{{n-1\choose i-1}(i-1)!  (d-1)^i d^{i-1}} 
{(nd-1)(nd-3)\cdots (nd-2i+1)}\\  \nonumber
&\leq& \frac{2d-4}{nd-1}+ \sum_{i=2}^{g-1} \frac{4 (d-1)^i}{nd-1}\\
&<& \epsilon. \label{Mj=1}
\end{eqnarray} 
Consider now a fixed  $M\in \M_j$  for some $j\geq 2$. We see $j$ classes and these $j$ classes and the edges of $M$ project to a $j$-cycle.
 If $M'\in \M_1$, then the single edge of $M'$ connects two vertices of the same class (one of the $j$ classes), such that one of its endpoints is an endpoint 
 of an edge of $M$ as well. There are $2d-3$ such edges in every class, totaling $j(2d-3)$. If $M'\in \M_i$ for some $i\geq 2$, then there is class $C$ containing two
 endpoints, $x$ and $y$, of two different edges of $M$, such that an edge of $M'$ has $x$ or $y$ as an endpoint. To build all such $M'$ matchings, select the $i-1$ classes in $\binom{n-1}{i-1}$ ways, and proceed as in the previous argument---there could be more overlapping with the  $j$ classes, but we only need an upper bound. 
 In conclusion, we obtain 
\begin{eqnarray}
    \sum_{M':A_{M'}\cap A_M=\emptyset}\Pr(A_{M'}) &=& \nonumber
\frac{j(2d-3)}{nd-1}+ 
\sum_{i=2}^{g-1} \sum_{M'\in \M_i: A_{M'}\cap A_M=\emptyset}\Pr(A_{M'})\\
&\leq& \frac{j(2d-3)}{nd-1} \nonumber
+ j\sum_{i=2}^{g-1} \frac{{n-1\choose i-1}(i-1)!  (d-1)^{i} d^{i-1}} 
{(nd-1)(nd-3)\cdots (nd-2i+1)}\\ \nonumber
&<& \frac{(2d-3)(g-1)}{nd-1}
+ (g-1)\sum_{i=2}^{g-1} \frac{ 4(d-1)^{i}} 
{nd-1}\\
&<& \epsilon. \label{Mjnot=1}
\end{eqnarray}
Now we verify item 4. Any single $uv$ edge can be in at most one $M\in \M_1$,  whose probability is $\frac{1}{nd-1}$. If $uv\in M\in \M_i$ for $i\geq 2$, then
in addition to the  classes of $u$ and $v$ we have select $i-2$ classes out of the remaining $n-2$. We can put the $i$ classes into a cycle such that the  classes of $u$ and $v$
are neighbors in $(i-2)!$ ways. Selecting the the endpoints of the $i-1$ edges different from $uv$ from the classes can be done in $(d-1)^i d^{i-2}$ ways. 
In conclusion, we obtain
\begin{eqnarray}
  \sum_{M: uv\in M\in\M}\Pr(A_{M})&\leq& 
\frac{1}{nd-1} \nonumber
+ \sum_{i=2}^{g-1} \frac{{n-2\choose i-2}(i-2)!  (d-1)^i d^{i-2}} 
{(nd-1)(nd-3)\cdots (nd-2i+1)}\\   \nonumber
&<& \frac{1}{nd-1}
+ \sum_{i=2}^{g-1}  \frac{1}{d(n-1)} \cdot \frac{ 4(d-1)^i} 
{nd-1}\\
&<&\epsilon/2. \label{containment}
\end{eqnarray}

Finally, we verify item 5.
For any $F\in \M$,  we need estimate \linebreak $\sum_{M\in \M_F} \Pr(A_M^{N-2r})$.
If the projection of $F$ is a loop, then $\M_F=\emptyset$ and there is nothing to do.
Now we assume the projection of $F$ is a cycle $C_k$. Assume that $M'\in \M$
intersects $F$, $M=M'\setminus F$,
and the projection of $M'$ is a cycle $C_s$ with $k,s\leq g-1$.
Then the components of $C_s\cap C_k$ having at least one edge are  paths $P_1, P_2,\ldots, P_t$, with  $t\geq 1$.
Fixing these paths, and the  edges in $M'\cap F$,
some additional $\ell$ vertices are joined with these $t$ paths
to make $C_s$. So the number of possible $C_s$'s 
with these fixed
paths 
  is at most
$$\sum_{\ell\leq g-1-2t} {n\choose \ell} (\ell+t-1)!2^t,$$
and the number of $M'$-s  defining this particular $C_s$   with  $M'\cap F$ fixed,
is at most 
$d^{\ell}(d-1)^{\ell+2t}$. The $t$ paths with at least one edge can be
selected  in at most $2 {k\choose 2t} $ ways from $C_k$. The probability
$\Pr(A_M^{N-2r})$, where $M= M'\setminus F$, is at most $(N-3g)^{-(\ell+t)}$.
 We summarize that 
\begin{equation} \label{tobound}
\sum_{M\in \M_F} \Pr(A_M^{N-2r})
\leq 
\sum_{t=1}^{\lfloor k/2  \rfloor }2 { k\choose 2t} \sum_{\ell\leq g-1-2t} {n\choose \ell} 
(\ell+t-1)!2^t \frac{d^{\ell}(d-1)^{\ell+2t}}{  (N-3g)^{\ell+t}}.
\end{equation}
As $\ell+t-1\leq g-3$, using the falling factorial notation we have $(\ell+t-1)!=\ell!(\ell+t-1)_{t-1}\leq \ell!(g-3)^{t-1}$.
There is an absolute upper bound $K>\frac{(n)_\ell d^\ell}{(N-3g)^\ell}$. As $\ell+2t\leq g-1$, the RHS of 
 (\ref{tobound}) can be further estimated by 
 \[2K(d-1)^{g-1}\sum_{t=1}^{\lfloor k/2 \rfloor}{k\choose 2t} \sum_{\ell\leq g-1-2t}\Biggl(\frac{2(g-3)   }{N-3g   }\Biggl)^t\leq   2Kg(d-1)^{g-1}\sum_{t=1}^{\lfloor k/2 \rfloor}{k\choose 2t}  \Biggl(\frac{2(g-3)   }{N-3g   }\Biggl)^t.   \]
It is easy to see that the last summation has the largest term at $t=1$, has less than 
$g$ terms,   
 and is $\leq 4Kg^5(d-1)^{g-1}/(N-3g)<\epsilon$.

To  apply Theorem \ref{enum}, we need $r\epsilon=o(1)$ and   $r\mu\epsilon=o(1)$. 
The first condition  follows from  the second as $\mu$ is separated from zero.  As $r<g$, $\mu\leq 
(d-1)^{g-1}/2$ and $\epsilon=\frac{K'g^5(d-1)^{g-2}}{n}$, the second condition boils down
to $g^6(d-1)^{2g-3}=o(n)$, which was provided in (\ref{condition}). The neglection
of error in (\ref{errorinmainterm}) is also allowed by  (\ref{condition}).
\hfill $\square$

In the {\em bipartite configuration model} we have two sets, $U$ and $V$, each 
containing $N$ mini-vertices,  a fixed partition of $U$ into $d_1,...,d_n$ element classes,
and a fixed partition of $V$ into $\delta_1,...,\delta_n$ element classes.
Any perfect matching between $U$ and $V$ defines a bipartite multigraph with
partite sets of size $n$ after a projection contracts every class to single vertex.
In the regular case, $d_1=\cdots=d_n=\delta_1=\cdots=\delta_n=d$.
We prove next another theorem of McKay, Wormald and Wysocka
\cite{wysocka}: 

\begin{theorem} \label{multibigraph}  In the regular case of the bipartite configuration model, assume that 
$g$ is even, $d\geq 3$, and 
\begin{equation} \label{condition2}
g^6(d-1)^{2g-3}=o(n).
\end{equation}
Then
the probability that the random bipartite $d$-regular multigraph has
girth at least $g\geq 2$ is 
$  (1+o(1))\exp\Bigl({-\sum_{i=1}^{(g-2)/2}\frac{(d-1)^{2i}}{2i}} \Bigl) $,
and hence the number of $d$-regular bipartite graphs on $n,n$ vertices with girth
 at least $g\geq 4$ is
\[(1+o(1))e^{-\sum_{i=1}^{(g-2)/2}\frac{(d-1)^{2i}}{2i}}\frac{(dn)!}{(d!)^{2n}}. \]
(The case $g=4$ means that the random $d$-regular bipartite multigraph is actually
a simple bipartite graph.) Furthermore, the number of $d$-regular bipartite graphs not containing 
cycles whose length is in a set ${\cal C}\subseteq\{4,6,...,g-2\}$, is
 \[(1+o(1))e^{ - \frac{(d-1)^2}{2}-\sum_{i\in {\cal C}}\frac{(d-1)^i}{i}}\frac{(dn)!}{(d!)^{2n}}
. \]
\end{theorem}
{\bf Proof:} We outline the proof of the first claim. 
For $i=1,2,\ldots, (g-2)/2$, let $\M_i$ be the set of  matchings of $U$ and $V$, whose
projection gives a cycle of length $2i$; there are {\em exactly}
${n\choose i}^2  d^{2i} (d-1)^{2i}(i-1)!^2i$
of them.  The bad events for the negative dependency
graph are the union of matchings $\M=\cup_{i=1}^{(g-2)/2} \M_i$.
For each $M_i\in \M_i$ ($i=1,2,\ldots, (g-2)/2$), we have
\begin{equation}
  \label{eq:mi2}
  \Pr(A_{M_i})=\frac{(dn-2i)!}{(dn)!}.
\end{equation}
We have
\begin{eqnarray}
&&\sum_{M\in \M}\Pr(A_M) =\sum_{i=1}^{(g-2)/2} {n\choose i}^2
 d^{2i} (d-1)^{2i}(i-1)!^2i
\nonumber
\frac{(dn-2i)!}{(dn)!}\\
&=&  \sum_{i=1}^{(g-2)/2}   \frac{(d-1)^{2i}}{2i}\left(1+ O\left(\frac{i^2}{n}\right)\right)
= \left(1+ O\left(\frac{g^2}{n}\right)\right)   \sum_{i=1}^{(g-2)/2} \frac{(d-1)^{2i}}{2i}   
\nonumber .\\ \label{errorinmainterm2}
\end{eqnarray}
All the estimates go through as in the proof of Theorem~\ref{multigraph}.  
 To prove the second claim, only
(\ref{errorinmainterm2}) has to be adjusted, everything else remains the same.
\hfill $\square$

\subsection{Applications III: Enumeration of graphs by girth and degree sequence}
McKay and Wormald \cite{81} enumerated graphs by degree sequences. We extend 
their result to include the girth or the set of allowed short cycle lengths. However,
our range for the degrees is not as broad as in \cite{81}. For example, formula 
(\ref{strongest}) that we could not obtain is a special case of \cite{81}.

We start with some technicalities on estimating elementary symmetric polynomials.
 Let $\sigma_n^{(k)}(x_1,...,x_n)$ denote the
$k^{th}$ elementary symmetric polynomial in $n$ variables. Assume that
every $x_i>0$ and set average $\bar x=(\sum_{i=1}^n x_i)/n$
and the second order average  ${\tilde x}=(\sum_{i=1}^n x_i^2)/\bar x$. We claim the following:
\begin{equation} \label{symmetric}
\frac{n^k}{(n)_k}\Bigl(1- \frac{{k \choose 2}}{n^2}\cdot \frac{{\tilde x}}{\bar x}\Bigl)
  \leq\frac{\sigma_n^{(k)}(x_1,...,x_n)}{\sigma_n^{(k)}(\bar x,...,\bar x)}\leq \frac{n^k}{(n)_k}.
\end{equation}
Now we verify  (\ref{symmetric}). First observe the inequality
 $$\sigma_n^{(k)}(x_1,...,x_n)\leq \frac{(x_1+...+x_n)^k}{k!},$$ which holds termwise for the two multivariate polynomials.
 This inequality immediately implies the upper bound in    (\ref{symmetric}). Next observe that
 $$    \frac{(x_1+...+x_n)^k-{k \choose 2}(\sum_{i=1}^n x_i^2)
(x_1+...+x_n)^{k-2}  }{k!}\leq         \sigma_n^{(k)}(x_1,...,x_n),               $$
 as  the inequality holds termwise for the two multivariate polynomials. This implies 
the lower bound in    (\ref{symmetric}). We conclude from   (\ref{symmetric})    the asymptotic formula
\begin{equation}\label{symmasymp}
\sigma_n^{(k)}(x_1,...,x_n)=\frac{n^k(\bar x)^k}{k!}\Biggl(1+O\Bigl(\frac{k^2}{n^2}\cdot
\frac{\tilde x}{\bar x}\Bigl)\Biggl),
\end{equation}
whenever the quantity in the $O$-term goes to zero. Assume further that $x_1\leq
x_2\leq ...\leq x_n$. Define a sequence by $y_i=x_{t+i}$ for $i=1,2,...,n-t$. It is easy 
to see the following chains of inequalities:
$$\bar x=\frac{x_1+...+x_n}{n}\leq \frac{x_{t+1}+...+x_n}{n-t}=\bar y \leq  \frac{x_{1}+...+x_n}{n-t}=\Bigl(1+\frac{t}{n-t}  \Bigl) \bar x$$
and
$$\frac{\frac{n-t}{n}(x_1^2+...+x_n^2)}{\frac{n}{n-t}\bar x}\leq \frac{x_{t+1}^2+...+x_n^2}{\bar y}=\tilde y\leq \frac{x_{1}^2+...+x_n^2}{\bar x}=\tilde x.$$
Based on them,  for $t=o(n)$  
we have
 $\bar y=\biggl(1+O\bigl( \frac{t}{n}\bigl)\biggl){\bar x}$
and $\tilde y=\biggl(1+O\bigl( \frac{t}{n}\bigl)\biggl){\tilde x}$. From here and
(\ref{symmasymp}) we conclude that $kt=o(n)$ implies 
\begin{equation}\label{symmasymp2}
\sigma_{n-t}^{(k)}(y_1,...,y_{n-t})=
\frac{n^k(\bar x)^k}{k!}\Biggl(1+O\Bigl(\frac{kt}{n}+\frac{k^2}{n^2}\cdot \frac{\tilde x}
{\bar x}\Bigl)\Biggl).
\end{equation}
To verify (\ref{symmasymp2}), observe
\begin{eqnarray*}
\sigma_{n-t}^{(k)}(y_1,...,y_{n-t})&=&  \frac{(n-t)^k(\bar y)^k}{k!} \Biggl(1+O\Bigl(\frac{k^2}{n^2}\cdot \frac{\tilde y}{\bar y}\Bigl)\Biggl)\\
&=& \frac{n^k(\bar x)^k}{k!} \Bigl( \frac{n-t}{t}  \Bigl)^k      \Bigl(1+O\bigl( \frac{t}{n}\bigl)  \Bigl)^k \Biggl(1+O\Bigl(\frac{k^2}{n^2}\cdot \frac{\tilde y}{\bar y}\Bigl)\Biggl).
\end{eqnarray*}

Let us return to the configuration model as described
at the beginning of Subsection~\ref{dreg} and try to do in more generality 
the steps of the proof of Theorem~\ref{multigraph}.  The combinatorial structures
are the same, but different class sizes have to be taken into account.
Assume now that $2\leq d_1\leq d_2\leq ... \leq d_n$ and set $D_j=d_j(d_j-1)$.
If the projection provides a {\em graph}
with degree sequence $d_1,d_2,...,d_n$ (as opposed to a multigraph), then
{\em exactly} $d_1!d_2!\cdots d_n!$ matchings on the set of $N=d_1+...+d_n$
mini-vertices yield this graph. We want to compute the probability that after the
projection we obtain a graph with girth at least $g$ ($g\geq 3$).
For $i=1,2,\ldots, g-1$, let $\M_i$ be the set of  matchings of $U$ whose
projection gives a cycle of length $i$; there are {\em exactly}
$$\frac{(i-1)!}{2} \sigma_n^{(i)}  \bigl(D_1 ,...,D_n\bigl)$$
of them. Assume $i\geq 3$. Consider an arbitrary $i$-cycle after the projection.
The $i$ vertices of a cycle must have come from $i$ disjoint classes.
Denote by $C_1,...,C_i$  those disjoint classes, in the cyclic  order of the vertices of the cycle.
Count how many matchings project to this fixed cycle.
Select a vertex in $C_1$ in $|C_1|$ ways, join it to vertex of $C_2$ in $|C_2|$ ways, select a second vertex of  $C_2$ in $|C_2|-1$ ways, join it to a 
vertex of $C_3$ in $|C_3|$ ways, select a second vertex of  $C_3$ in $|C_3|-1$ ways, ... , etc., ... , select a second vertex of  $C_i$ in $|C_i|-1$ ways, join it to a 
vertex of $C_1$ in $|C_1|-1$ ways. We found $\prod_{j=1}^i \Bigl( |C_j|\cdot (  |C_j|-1)  \Bigl)$ ways. This number of ways is independent of the cyclic order of classes.
Using all $(i-1)!$ cyclic orders, however, we obtain every cycle exactly twice---going to the left from $C_1$ and going to the right from $C_1$. 
It is nice and easy to verify that for the  degenerate cases $i=1,2$ the same formula works.
 
 The bad events for the negative dependency
graph are the union of matchings $\M=\cup_{i=1}^{g-1} \M_i$.
For each $M_i\in \M_i$ ($i=1,2,\ldots, g-1$), we have 
\begin{equation}
  \label{eq:mi3}
 \Pr(A_{M_i})=\frac{1}{(N-1)(N-3)\cdots (N-2i+1)},
\end{equation}
where $N=n\bar d$. We find under the mild assumption $  \frac{g^2}{n}+
\frac{g^2}{n^2}\cdot \frac{\tilde D} 
{\bar D}=o(1)     $ that 
\begin{eqnarray}
&&\sum_{M\in \M}\Pr(A_M) = \sum_{i=1}^{g-1} \frac{(i-1)!}{2}\cdot
\nonumber
\frac{ \sigma_n^{(i)}  \bigl(D_1,...,D_n\bigl)   }{ (N-1)(N-3)\cdots (N-2i+1)  }\\
\hbox{(by (\ref{symmasymp}))}&=&     \sum_{i=1}^{g-1} 
\frac{ n^i  (\bar D)^i  }{ 2i(N-1)(N-3)\cdots (N-2i+1)  }\Biggl(1+O\Bigl(
\frac{i^2}{n^2}\frac{\tilde D}
{\bar D}\Bigl)\Biggl) 
\nonumber \\ 
&=&     \sum_{i=1}^{g-1} 
 \frac{1}{2i} {\Bigl(\frac{\bar D}{\bar d}\Bigl)^i  } \Biggl(1+O\Bigl(\frac{i^2}{n}+
\frac{i^2}{n^2}\frac{\tilde D}
{\bar D}\Bigl)\Biggl)   \nonumber \\
&=&\Biggl(1+O\Bigl(\frac{g^2}{n}+
\frac{g^2}{n^2}\frac{\tilde D} 
{\bar D}\Bigl)\Biggl)\sum_{i=1}^{g-1} 
 \frac{1}{2i} {\Bigl(\frac{\bar D}{\bar d}\Bigl)^i  }. 
\label{errorinmaintermgen}
\end{eqnarray}
The estimate in (\ref{Mj=1}) changes to
\begin{equation} \label{Mj=1gen}
\frac{2d_n-4}{n{\bar d}-1}+ \sum_{i=2}^{g-1} \frac{(i-1)!  (d_n-1)
\sigma_{n-1}^{(i-1)}(D_2,...,D_n)} 
{(n{\bar d}-1)(n{\bar d}-3)\cdots (n{\bar d}-2i+1)}. 
\end{equation}
To see this, go back to the proof of  (\ref{Mj=1}). There is a fixed $\{e\}=M\in \M_1$ with $e=xy$ belongs to the same class $C_1$ that
projects to a mini-vertex. $M'$ is a matching of $i$ edges that projects to an $i$-cycle, and $M'$ has at least one of $x$ and $y$ among the endpoints of its edges.
Assume that $C_1,...,C_i$ are the classes that are the pre-images of the vertices of the $i$-cycle, in this cyclic order. Following
 the cyclic order, exactly $|C_2|\cdot (|C_2|-1)\cdots |C_i|\cdot (|C_i|-1) (|C_1|-1)$ matchings project to this cycle.  Every cycle, however, can obtained twice from this
 procedure, from two mirror image cyclic orders. At this point we estimate by $|C_i|\leq d_n$.  The classes $C_2,...,C_i$ are selected from the remaining  $n-1$ classes.
 If $C_1$ was the class with index $j$ and $d_j$ vertices from the list of all classes, the sum of $|C_2|\cdot (|C_2|-1)\cdots |C_i|\cdot (|C_i|-1)$ for all selections of $i-1$ 
 classes is exactly
 $$\sigma_{n-1}^{(i-1)}(D_1,..., \widehat{D_j},...,D_n) \leq
\sigma_{n-1}^{(i-1)}(D_2,...,D_n)$$
as the $D$ sequence is increasing.
The estimate in (\ref{Mjnot=1}) changes to
\begin{equation} \label{Mjnot=1gen}
\frac{j(2d_n-3)}{n{\bar d}-1} 
+ j\sum_{i=2}^{g-1} \frac{(i-1)!  (d_n-1) \sigma_{n-1}^{(i-1)}(D_2,...,D_n) } 
{(n{\bar d}-1)(n{\bar d}-3)\cdots (n{\bar d}-2i+1)}
\end{equation}
by an almost identical argument, following the proof of (\ref{Mjnot=1}).
The estimate in (\ref{containment}) will change to
\begin{equation} \label{containmentgen}
\frac{1}{n{\bar d}-1} 
+ \sum_{i=2}^{g-1} \frac{(d_n-1)^2 (i-2)!\sigma_{n-2}^{(i-2)}(D_3,...,D_n) } 
{(n{\bar d}-1)(n{\bar d}-3)\cdots (n{\bar d}-2i+1)}.
\end{equation}
Indeed, mimicking the proof of (\ref{containment}), assume that an $uv$ edge connects classes $C_1$ and $C_i$, where 
$C_1$ was the class with index $p$ and $d_p$ vertices,  while $C_i$ was the class with index $q$ and $d_q$ vertices  from the list of all classes, $p\not= q$.
Fix an $i$-cycle after the projection that contains the mini-vertices arising from $C_1$ and $C_i$. Assume the classes corresponding to the mini-vertices in the order
of the cycle are $C_1,C_2,...,C_i$.
Any $i$-matching projecting to the fixed $i$-cycle can come into existence $(|C_1|-1) \cdot |C_2|\cdot (|C_2|-1)\cdots |C_{i-1}|\cdot (|C_{i-1}|-1)\cdot (|C_i|-1)$ ways.
We estimate our terms by $(|C_1|-1)(|C_i|-1)\leq (d_n-1)^2$ and 
 $$\sigma_{n-2}^{(i-2)}(D_1,..., \widehat{D_p},...,  \widehat{D_q}  ,...,D_n) \leq
\sigma_{n-2}^{(i-2)}(D_3,...,D_n).$$

The estimate in (\ref{tobound}) changes to
\begin{equation} \label{toboundgen}
\sum_{t=1}^{\lfloor k/2  \rfloor }2 { k\choose 2t} \sum_{\ell\leq g-1-2t} 
 \frac{ (\ell+t-1)!2^t    (d_n-1)^{2t}\sigma_{n-2t}^{(\ell)} (D_{2t+1},...,D_n) }{  (N-3g)^{\ell+t}},
 \end{equation}
 with an explanation for the elementary symmetric polynomial like at the last three numbered formulae.
 We are in a position to claim to the generalization of Theorem~\ref{multigraph} for other than constant degree sequences. It is remarkable that we do not have to assume the
 Erd\H os-Gallai condition for the targeted sequence, as our conditions imply it.
 \begin{theorem} \label{degmultigraph} Assume that $N=d_1+...+d_n$ is even,
 $\bar d \geq 3$, every $d_i\geq 2$.
  In the configuration model, assume 
\begin{equation} \label{condition3}
\Bigl( \frac{g^2}{n}+ \frac{g^2}{n^2}\cdot \frac{\tilde D}{\bar D}\Bigl)\cdot \Bigl(\frac{\bar D}{\bar d}\Bigl)^{g-1}=o(1) \hbox{\ \ and \ \ }
g^6   \Bigl(  \frac{\bar D}{\bar d} \Bigl)^{2g-4}d_n^2=o(N).
\end{equation}
Then
the probability that the random  multigraph  with degrees $d_1,d_2,...,d_n$ after the 
projection  has
girth at least $g\geq 1$ is 
\begin{equation}  \label{girthprob}
(1+o(1))\exp\Biggl(-\sum_{i=1}^{g-1}\frac{1}{2i}\cdot 
 \biggl(\frac{\bar D}{\bar d}\biggl)^i  
  \Biggl),
\end{equation}
and hence the number of  graphs on $n$ vertices with  degrees   $d_1,d_2,...,d_n$ and  girth
 at least $g\geq 3$ is
\[(1+o(1)) \frac{(N-1)!!}{ \prod_i d_i!}
\exp\Biggl(-\sum_{i=1}^{g-1}\frac{1}{2i}\cdot 
 \biggl(\frac{\bar D}{\bar d}\biggl)^i  
  \Biggl)
   .\]
(The case $g=3$ means that the random  multigraph is actually
a simple graph, and  hence  $d_1,d_2,...,d_n$ is a graph degree sequence.) Furthermore, the number of  graphs    with  degrees   $d_1,d_2,...,d_n$    not containing 
cycles whose length is in a set ${\cal C}\subseteq\{3,4,...,g-1\}$, is
 \[(1+o(1))\frac{(N-1)!!}{\prod_id_i!}
\exp\Biggl( -\frac{1}{2}\cdot 
 \biggl(\frac{\bar D}{\bar d}\biggl)
 -\frac{1}{4}\cdot 
 \biggl(\frac{\bar D}{\bar d}\biggl)^2
 -\sum_{i\in {\cal C}}
 \frac{1}{2i}\cdot 
 \biggl(\frac{\bar D}{\bar d}\biggl)^i  
\Biggl)
. \]
\end{theorem}
{\bf Proof:} Take  $\epsilon=
Kg^5\bigl(\frac{\bar D}{\bar d}\bigl)^{g-3}\frac{(d_n-1)^2}{N-3g}$
for some large constant  $K$.
The estimate to (\ref{toboundgen}) goes similar to the estimate for (\ref{tobound}), but to estimate the elementary symmetric polynomial it uses 
(\ref{symmasymp2}):
$$\leq\sum_{t=1}^{\lfloor k/2  \rfloor }2 { k\choose 2t} \sum_{\ell\leq g-1-2t} 
 \Biggl[    \frac{4(g-3)(d_n-1)^2}{N-3g}           \Biggl]^t \Biggl(\frac{\bar D}{\bar d}\Biggl)^\ell
 \Bigl(1+O\bigl(\frac{\ell t}{n}+  \frac{\ell^2}{n^2}\cdot \frac{\bar D}{\tilde D} \bigl)\Bigl).
$$
The second part of (\ref{condition3}) implies $g^3d_n^2=O(N)$, which in turn implies that
the biggest term in the bound occurs for $t=1$. There are at most $g$ terms, and therefore
$\epsilon$ bounds (\ref{toboundgen}).
We leave it to the reader that this $\epsilon $ also provides a bound for (\ref{containmentgen}),
 (\ref{Mjnot=1gen}), and  (\ref{Mj=1gen}). The least trivial is the middle one, it follows from inequality 
  $\frac{\bar D}{\bar d}\leq d_n$.

The Cauchy-Schwartz inequality shows that 
\begin{equation}\label{cauchy}
\sum_{i=1}^{g-1}\frac{1}{2i}\cdot 
 \bigl(\frac{\bar D}{\bar d}\bigl)^i=O\Bigl(  \bigl(\frac{\bar D}{\bar d}\bigl)^{g-1} \Bigl)
 \end{equation}
and  therefore the first part of (\ref{condition3}) allows the approximation
 in (\ref{errorinmaintermgen}). The second part of (\ref{condition3}) implies  the conditions
 above (\ref{eq:enum}) and  (\ref{eq:enum2}).

 The proof of the second
claim is analogous.
 \hfill $\square$

It is not difficult to obtain a degree sequence version of Theorem~\ref{multibigraph}.
As the proof is just a combination of the proofs of Theorems~\ref{multibigraph} 
and~\ref{degmultigraph}, we leave the details to the reader.

\begin{theorem} \label{degmultibigraph}  In 
the bipartite configuration model, assume that 
$g$ is even, the class sizes are $2\leq d_1\leq\cdots \leq d_n$ and $
2\leq \delta_1\leq \cdots \leq \delta_n$,
$N=\sum_i d_i=\sum_i \delta_i$,
${\bar d}={\bar \delta}\geq 3$, $D_j=d_j(d_j-1)$ and $\Delta_j=\delta_j(\delta_j-1)$.
Assume further that 
\begin{equation} \label{condition5}
\frac{g^2}{n^2}\Bigl(n+ \frac{\tilde D}{\bar D} +  \frac{\tilde \Delta}{\bar \delta} \Bigl)  
\Bigl(\frac{\bar D\cdot \bar \Delta}{\bar d\cdot \bar\delta}  \Bigl)^{(g-2)/2}=o(1) \hbox{\ \ and \ \ } g^6(d_n^2+\delta_n^2)
\Bigl( \frac{\bar D}{\bar d} \Bigl)^{g-3}\Bigl( \frac{\bar \Delta}{\bar \delta} \Bigl)^{g-3}
=o(N).
\end{equation}
Then
the probability that the random bipartite  multigraph with the prescribed degree
sequence has
girth at least $g\geq 2$ is 
$$  (1+o(1))\exp\Bigl(-\sum_{i=1}^{(g-2)/2}\frac{(\bar D)^i(\bar \Delta)^i}{2i(\bar d)^{2i}} \Bigl), $$
and hence the number of bipartite graphs   with the prescribed degree
sequence and girth
 at least $g\geq 4$ is
\[(1+o(1))\frac{N!}{\prod_i d_i!\delta_i!}\exp\Bigl(-\sum_{i=1}^{(g-2)/2}
\frac{(\bar D)^i(\bar \Delta)^i}{2i(\bar d)^{2i}}
\Bigl). \]
(The case $g=4$ means that the random  bipartite multigraph with the  given
degree sequence is actually
a simple bipartite graph, and hence given sequence {\em is} a bipartite graph degree 
sequence.) 
Furthermore, the number of  bipartite graphs with the prescribed degree 
sequence that do  not contain
cycles whose length is in a set ${\cal C}\subseteq\{4,6,...,g-2\}$, is
 \[(1+o(1))\frac{N!}{\prod_i d_i!\delta_i!}
\exp\Bigl(
 - \frac{{\bar D}{\bar \Delta}}{2(\bar d)^2}-\sum_{i\in {\cal C}}
 \frac{(\bar D)^i(\bar \Delta)^i}{2i(\bar d)^{2i}}
  \Bigl). \]
\end{theorem}
\section{Revisiting girth and chromatic number: high girth and high chromatic number graphs on a 
 given degree sequence}
An early result of Erd\H{o}s \cite{erdoscan}
asserts that for every $k$ and $g$, there is a graph $G$ with
$girth(G)\geq g$ and chromatic number $\chi(G)\geq k$. 
In Theorem~\ref{ujerdos} we refine this result of Erd\H os, changing the existential
quantifier to universal.

We start with some technicalities.
 Let $N$ be an even positive integer.
For a set $S\subset [N]$, we say that a perfect matching $M$ of $K_N$ 
{\em traverses} $S$, if every edge in $M$ is incident to at most 
one vertex in $S$, in other words no edge has two endpoints in $S$.
\begin{lemma}\label{transversal}
For a fixed set $S$ of size $s$, the probability that
$S$ is traversed by a random matching, equals to
$$\frac{2^s {{\frac{N}{2} \choose s}}}{{N \choose s}}.$$
This number is less than $e^{-\frac{s(s-1)}{2N}}$.
\end{lemma}
{\bf Proof:} Clearly the probability in question does not depend on the
choice of $S$, just depends on the cardinality $s$. Therefore the 
probability does not change if we average it out for all $s$-subsets,
and hence it is 
\[ 
\frac{\#(S,M): \hbox{\rm perfect matching $M$ traverses $S$}}
{(N-1)!!{N\choose s}}.
\]
Count now in the ordered pairs in the numerator as follows: for all
$(N-1)!!$ perfect matchings, decide which $s$ edges of the $N/2$ edges
of the perfect matching have endpoint in $S$, and for those $s$ edges 
decide which endpoint out of the two possibilities will belong to $S$.
For the estimate, \\
$2^s\frac{(N/2)_s}{(N)_s}=\prod_{i=0}^{s-1} \frac{N-2i}{N-i}
\leq\exp\Bigl(-\sum_{i=0}^{s-1} \frac{i}{N-i}\Bigl)\leq \exp\Bigl(-\frac{1}{N}\sum_{i=0}^{s-1} {i}\Bigl)= e^{-\frac{s(s-1)}{2N}}$.
\hfill $\square$

\begin{theorem} \label{ujerdos}
Consider the configuration model as in Theorem~\ref{degmultigraph}.
Assume (\ref{condition3}), $\bar d\geq 3$, fix $k$ and $\epsilon >0$, such that
 $k^2<(1-\epsilon)\frac{\bar d }{2\log 2}$.
Then almost all graphs with degree sequence $d_1,...,d_n$ and girth at least $g\geq 3$
are not $k$-colourable.
\end{theorem}
Specializing to regular graphs, we get back the existence of graphs of high 
chromatic number and high girth, roughly in the same range where Erd\H os 
\cite{erdoscan} obtained it. 

\noindent{\bf Proof.}  Consider a random matching $M$ on $N$ vertices and its contraction 
into a multigraph $G$ with the prescribed degree sequence. We have to show
$$\Pr\bigl(\chi(G)\leq k \bigl\vert G \hbox{\ simple}\bigl)=o\biggl(\Pr\bigl(girth(G)\geq g \bigl\vert 
G \hbox{\  simple}\bigl)\biggl).$$
This is equivalent to
\begin{equation} \label{whatweneed}
\Pr\bigl(G \hbox{\  is simple and }\chi(G)\leq k \bigl)=o\biggl(\Pr\bigl(girth(G)\geq g \bigl)\biggl).  
\end{equation}
Recall that (\ref{girthprob}) gave the probability that
the multigraph resulting from the configuration model has girth at least $g$.
Because of the $g\geq 3$ assumption, the probability that a resulting {\em graph}
has girth at least $g$ is  the same  (\ref{girthprob}).

Now we set an  upper bound on the probability that a simple 
$G$ is $k$-colorable.  For a subset  $A$ of $V(G)$, let the {\em volume} of $A$ 
be $\sum_{v\in A} d_G(v)$.
If $G$ is simple and $k$-colorable, then
$G$ contains an independent set of volume at least $\frac{N}{k}$.
By Lemma \ref{transversal},  at $s=\lceil N/k\rceil$, 
the probability of this event is at most
\begin{equation} \label{travers}
2^n \exp\Biggl(-\frac{N}{2k^2}+\frac{1}{2k}\Biggl)= \exp
\Biggl(\Bigl(-\frac{1}{2k^2}-\frac{\log 2}{ \bar d}
\Bigl) N+\frac{1}{2k}\Biggl ).
\end{equation}
Computing the difference of the exponents in (\ref{travers}) and in (\ref{girthprob})
we are at home, if we use (\ref{cauchy}) to bound the exponent in (\ref{girthprob}),
and the second part of (\ref{condition3}) with $2g-4\geq g-1$.
\hfill $\square$\\
\noindent{\bf Acknowledgement.} We thank  Austin Mohr and an anonymous referee  for their valuable comments.


\begin{thebibliography}{9}
\bibitem{as} N.~Alon,  J.~H.~Spencer, {\em The Probabilistic Method,}
Second Edition,
   John Wiley and Sons, New York, 2000.








\bibitem{bender} E.A.~Bender, Asymptotic methods in enumeration,
{\em SIAM Review}  {\bf 18}(4), (1974),  485--515 (see in particular pages 510--511).


\bibitem{benderca} E. A. Bender, E. R. Canfield, The asymptotic number of
non-negative integer matrices with given row and column sum, {\em 
J. Comb. Theory} Ser. A {\bf 24} (1978), 296--307.

\bibitem{configmodel} B.~Bollob\'as, A probabilistic proof of an 
asymptotic formula for the number of labelled regular graphs,
{\em Europ. J. Combinatorics} {\bf 1} (1980), 311--316.

\bibitem{boppanas} R. B. Boppana, J. H. Spencer, A useful elementary 
correlation inequality, {\em J. Comb. Theory} A {\bf 50} (1989), 305--307.

\bibitem{borkowitz} D.~Borkowitz, The name of the 
game: exploring random permutations,
{\em Mathematics Teacher} {\bf 98} (2005) (October) 
and its appendix\\
{\tt http://faculty.wheelock.edu/dborkovitz/articles/ngm6.htm}

\bibitem{chen} L.~Chen, Poisson approximation for independent trials,
{\em Ann. Probab.} {\bf 3} (1975), 534--545.



 \bibitem{erdoscan} P.~Erd\H{o}s, Graph theory and probability,
{\em  Canad. J. Math.} {\bf 11} (1959), 34--38.

\bibitem{erdoskap} P.~Erd\H{o}s, J.~Kaplansky, 
The asymptotic number of Latin rectangles, {\em Amer. J. Math.} {\bf 68} 
(1946), 230--236.

 \bibitem{erdoslovasz}
 P.~Erd\H{o}s and L.~Lov\'asz, Problems and results on 3-chromatic
 hypergraphs and some related questions,
 in {\em Infinite and Finite Sets}, A.~Hajnal et.~al., Eds.,
 {\em Colloq.~Math.~Soc. J. Bolyai} {\bf 11}, North Holland, Amsterdam,
 1975, 609--627.

\bibitem{erdosspencer} P.~Erd\H{o}s and J.H.~Spencer,
Lopsided Lov\'asz local lemma and latin transversals,
{\em Discrete Appl. Math.} {\bf 30} (1991), 151--154.

\bibitem{godsil} C.D.~Godsil, B. D.~McKay, Asymptotic enumeration
of Latin rectangles, {\em J. Comb. Theory} B {\bf 48} (1990), 19--44.


\bibitem{ku} C.Y.~Ku, Lov\'asz local lemma, 
{\tt http://www.maths.qmul.ac.uk/\~{}cyk/}

\bibitem{soda} L.~Lu and L.A.~Sz\'ekely,
Using Lov\'asz Local Lemma in the space of random injections, 
{\em Electronic J. Comb.} {\bf 14} (2007), R63.

\bibitem{JCMCC} Linyuan Lu, A. Mohr, L. A. Sz\'ekely, Connected balanced subgraphs in random regular multigraphs under the configuration model,  to appear in {\em J. Comb. Math.
Comb. Comput.} (2013)


\bibitem{72} B.D.~McKay, Asymptotics for symmetric 0-1 matrices
with prescribed row sums, {\em Ars Combinatoria} {\bf 19A} (1985), 15--25.


\bibitem{81} B.D.~McKay and N.C.~Wormald, Asymptotic enumeration
by degree sequence of graphs with degrees $o(n^{1/2})$, {\em Combinatorica}
{\bf 11} (1991), 369--382.

\bibitem{wysocka} B.D.~McKay, N.C.~Wormald, and B.~Wysocka,
Short cycles in random regular graphs,
{\em Electronic J. Combinatorics} {\bf 11} (2004), \#R66, 12 pages.

\bibitem{procacci}
R.~Bissacot, R.~Fern\'andez, A.~Procacci, and B.~Scoppola,
An improvement of the Lov«asz Local Lemma via cluster
expansion,
arXiv:0910.1824v2 [math.CO] 26 Mar 2010

\bibitem{riordan} J.~Riordan, Three-line rectangles, {\em Amer. Math.
Monthly} {\bf 68} (1946), 230--236.


\bibitem{scott}
A.D.~Scott and A.D.~Sokal,
The repulsive lattice gas, the independent-set polynomial, and the Lov\'asz Local Lemma,
{\em J. Stat. Physics} {\bf 118}(5--6), (2005), 1151--1261.

\bibitem{spencerbook}
J.H.~Spencer, {\em Ten Lectures on the Probabilistic Method},
CBMS {\bf 52}, SIAM, 1987.

\bibitem{skau} I.~Skau, A note on the asymptotic number of Latin
rectangles, {\em Europ. J. Combinatorics} {\bf 19} (1998), 617-620.



\bibitem{stein} C.~Stein, Asymptotic evaluation 
of the number of Latin rectangles, {\em J. Comb. Theory} A
 {\bf 25} (1978), 38--49.

\bibitem{yamamoto} K.~Yamamoto, On the asymptotic number of Latin rectangles,
{\em Japan  J. Math.}  {\bf 21} (1951), 113--119.

\bibitem{wilf} H.S.~Wilf, {\em Generatingfunctionology}, Academic Press, 1990 (see pages
147--148). 

\bibitem{wormaldos} N.C.~Wormald, The asymptotic distribution of short 
cycles in  random regular graphs,  {\em J. Comb. Theory} Series B {\bf 31}
(1981), 168--182.




\bibitem{wormald} N.C.~Wormald,  Models of random regular graphs. 
{\em Surveys in combinatorics,} 1999 (Canterbury),  239--298, 
London Math. Soc. Lecture Note Ser., {\bf 267}, Cambridge Univ. Press, 
Cambridge, 1999.

\end{thebibliography}
\end{document}